\documentclass[11pt,a4paper]{article}
\usepackage{amsmath}
\usepackage[a4paper]{geometry}
\usepackage{amssymb,latexsym,amsmath,amsfonts,amsthm}
\usepackage{graphicx}
\usepackage{float}
\usepackage{subfigure}
\usepackage{algorithm}
\usepackage{algorithmic}
\usepackage{dsfont}
\usepackage{mathrsfs}
\usepackage{epsfig}
\usepackage{epstopdf}
\usepackage{booktabs}
\usepackage{comment,verbatim}
\usepackage{overpic}
\usepackage{hyperref}
\usepackage{bbm}
\usepackage[toc,page]{appendix}
\usepackage{multirow,diagbox}

\newcommand{\ud}{\,\mathrm{d}}
\newtheorem{theorem}{Theorem}[section]

\newtheorem{corollary}[theorem]{Corollary}

\theoremstyle{definition}

\theoremstyle{definition}

\theoremstyle{remark}

\newtheorem{remark}[theorem]{Remark}

\usepackage{color}

\usepackage[normalem]{ulem}

\def\Xint#1{\mathchoice
{\XXint\displaystyle\textstyle{#1}}%
{\XXint\textstyle\scriptstyle{#1}}%
{\XXint\scriptstyle
\scriptscriptstyle{#1}}%
{\XXint\scriptscriptstyle
\scriptscriptstyle{#1}}%
\!\int}
\def\XXint#1#2#3{{
\setbox0=\hbox{$#1{#2#3}{\int}$}
\vcenter{\hbox{$#2#3$}}\kern-.5\wd0}}
\def\dashint{\Xint-}

\numberwithin{equation}{section}

\begin{document}
\title{Complex generalized Gauss-Radau quadrature rules for Hankel transforms of integer order}
\author{Haiyong Wang\footnotemark[1]~\footnotemark[2] ~~ and ~~  Menghan Wu\footnotemark[1] }
\date{}

\footnotetext[1]{School of Mathematics and Statistics, Huazhong University of Science and Technology, Wuhan 430074, P. R. China. \texttt{Email:haiyongwang@hust.edu.cn}}

\footnotetext[2]{Hubei Key Laboratory of Engineering Modeling and Scientific Computing, Huazhong University of Science and Technology, Wuhan 430074, P. R. China}

\maketitle

\begin{abstract}
Complex Gaussian quadrature rules for oscillatory integral transforms have the advantage that they can achieve optimal asymptotic order. However, their existence for Hankel transform can only be guaranteed when the order of the transform belongs to $[0,1/2]$. In this paper we introduce a new family of Gaussian quadrature rules for Hankel transforms of integer order. We show that, if adding certain value and derivative information at the left endpoint, then complex generalized Gauss-Radau quadrature rules that guarantee existence can be constructed and their nodes and weights can be calculated from a half-size Gaussian quadrature rule with respect to the generalized Prudnikov weight function. Orthogonal polynomials that are closely related to such quadrature rules are investigated and their existence for even degrees is proved. Numerical experiments are presented to show the performance of the proposed rules.
\end{abstract}


\vspace{0.05in}


\section{Introduction}\label{sec:introduction}
Hankel transform of the form
\begin{equation}\label{eq:HankelTrans}
(\mathcal{H}_{\nu}f)(\omega) := \int_{0}^{\infty} f(x) J_{\nu}(\omega x) \ud x,
\end{equation}
where $J_{\nu}(x)$ is the Bessel function of the first kind of order $\nu$ and $\omega$ is the frequency of oscillations, appears in many physical problems, e.g., the propagation of optical, acoustic and electromagnetic fields \cite{Christopher1991,Sherer2004}, electromagnetic geophysics \cite{Denich2023,Ghosh1971,Ingeman2006,Ward1987} and layered media Green's functions \cite{Wang2021}.
Closed form of this transform is rarely available and numerical methods are generally required. However, conventional numerical methods for \eqref{eq:HankelTrans} are quite expensive due to the oscillatory and possibly slowly decaying behaviors of the integrand.

The evaluation of the Hankel transform has received considerable attention due to its practical importance and several methods have been developed, such as the digital linear filter (DLF) method \cite{Anderson1975,Ghosh1971}, the integration, summation and extrapolation (ISE) method \cite{Longman1956,Lucas1995}, complex Gaussian quadrature rules \cite{Asheim2013,Wong1982}. Among these methods, the DLF and ISE methods have their own disadvantages that the former requires samples on an exponential grid and the latter requires the evaluation of the zeros of $J_{\nu}(x)$. Complex Gaussian quadrature rules are of particular interest due to the appealing advantage that their accuracy improves rapidly as $\omega$ increases. Specifically, Wong in \cite{Wong1982} considered the construction of complex Gaussian quadrature rules for Bessel transforms of the form
\[
\int_{0}^{\infty} x^{\mu} f(x) H_{\nu}^{(s)}(x) \ud x, \quad s=1,2,
\]
where $H_{\nu}^{(s)}(x)$ are the Hankel functions and $\mu\pm\nu>-1$. By rotating the integration path to the imaginary axis such that the integrand is non-oscillatory and decays exponentially, Wong constructed some complex Gaussian quadrature rules for the above transforms. However, when extending Wong's method to the transform \eqref{eq:HankelTrans}, the resulting integrand will involve a nonintegrable singularity at the origin for $\nu\geq1$. More recently, Asheim and Huybrechs in \cite{Asheim2013} studied the construction of complex Gaussian quadrature rule for the Hankel transform \eqref{eq:HankelTrans} as well as Fourier and Airy transforms. Their key idea is to construct Gaussian quadrature rules with respect to the oscillatory weight function directly. Taking the Hankel transform \eqref{eq:HankelTrans} for example, they introduced the sequence of monic polynomials $\{P_n\}_{n=0}^{\infty}$ that are orthogonal with respect to $J_{\nu}(x)$ on $(0,\infty)$, i.e.,
\begin{equation}\label{eq:PolyBes}
\int_{0}^{\infty} P_n(x) x^k J_{\nu}(x) \mathrm{d}x = 0 , \quad
k=0,\ldots,n-1,
\end{equation}
and the above improper integrals are defined in terms of their Abel limits (see subsection \ref{sec:Abel}). Since $J_{\nu}(x)$ is a sign-changing function, the existence and uniqueness of the polynomials $P_n(x)$ cannot be guaranteed. Once the polynomials exist, Gaussian quadrature rule for the Hankel transform \eqref{eq:HankelTrans} can be derived immediately by using a simple scaling. A remarkable advantage of such rules is that they have the optimal asymptotic order, i.e., their error decays at the fastest algebraic rate as $\omega\rightarrow\infty$ among all quadrature rules using the same number of points. In the special case $\nu=0$, it was shown in \cite[Theorem 3.5]{Asheim2013} that the zeros of $P_n(x)$ with even $n$ are located on the imaginary axis. For $\nu>0$, based on numerical calculations for $\nu=1/2,1,3/2,2$, it was observed that the zeros of $P_n(x)$ tend to cluster near the vertical line $\Re(z)=\nu\pi/2$ as $n\rightarrow\infty$ (see Figure \ref{fig:GRNodesInt} for more detailed observations). More recently, this observation was proved in \cite{Deano2016} for $\nu\in[0,1/2]$ by using the steepest descent method for the Riemann-Hilbert problem of $P_n(x)$. However, the existence of the polynomials $\{P_n\}_{n=0}^{\infty}$ for $\nu>1/2$ still remains an issue.

On the other hand, let $\mu+\nu>-1$ and let $\{P_n^{(\mu,\nu)}\}_{n=0}^{\infty}$ be the sequence of monic polynomials that are orthogonal with respect to the weight $x^{\mu}J_{\nu}(x)$ on $(0,\infty)$, i.e.,
\begin{equation}\label{eq:PolyBesGen}
\int_{0}^{\infty} P_n^{(\mu,\nu)}(x) x^{k+\mu} J_{\nu}(x) \ud x = 0, \quad  k=0,\ldots,n-1.
\end{equation}
It is easily seen that such polynomials include the polynomials in \eqref{eq:PolyBes} as a special case and their existence cannot be guaranteed since the weight $x^{\mu}J_{\nu}(x)$ is also a sign-changing function. Unexpectedly, we prove that, if $\mu-\nu$ is a nonnegative integer, such polynomials exist for all degrees when $\mu-\nu$ is even and for only even degrees when $\mu-\nu$ is odd. Moreover, in these cases where the polynomials exist, they can be expressed in terms of the generalized Prudnikov polynomials. In this work, we consider a new family of Gaussian quadrature rules with theoretical guarantees for Hankel transform \eqref{eq:HankelTrans} of integer order. We begin by giving conditions for rotating the integration path of oscillatory integral transforms, including Hankel and Fourier transforms, in the right-half plane under their Abel limits. We then show that, if adding certain value and derivative information of $f(x)$ at the left endpoint, complex Gaussian quadrature rules for the Hankel transform \eqref{eq:HankelTrans} of integer order can be constructed with guaranteed existence. Since such rules involve the value and derivative information of $f(x)$ at the left endpoint and their nodes all lie on the imaginary axis, we refer to them as complex generalized Gauss-Radau quadrature rules. These theoretical findings on the polynomials defined in \eqref{eq:PolyBesGen} not only provide a theoretical rationale but are also useful for developing algorithms for the proposed complex generalized Gauss-Radau quadrature rules.

The rest of this paper is organized as follows. In the next section, we review some basic properties of generalized Gauss-Radau quadrature rules and the regularization of improper Riemann integrals. In Section \ref{sec:HankelInteger}, we construct complex generalized Gauss-Radau quadrature rules for the Hankel transform of integer order and study the existence of orthogonal polynomials that are closely related to such quadrature rules. In Section \ref{sec:Appl}, we give applications to oscillatory Hilbert transform and electromagnetic field configurations. Finally, we give some conclusions in Section \ref{sec:concluding}.

\section{Preliminaries}
We introduce some basics of generalized Gauss-Radau quadrature rules and the regularization of improper Riemann integrals. Throughout the paper, we denote by $\mathcal{P}_n$ the space of polynomials of degree at most $n$, i.e., $\mathcal{P}_n=\mathrm{span}\{1,x,\ldots,x^n\}$, and by $\mathbb{N}_0$ the set of nonnegative integers. Moreover, we denote by $\mathcal{K}$ a generic positive constant.

\subsection{Generalized Gauss-Radau quadrature rule}\label{sec:GenGaussRadau}
In this subsection we review some basic properties of generalized Gauss-Radau quadrature rule. Consider the following integral
\begin{equation}\label{eq:Int}
I(f) = \int_{a}^{b} w(x) f(x) \ud x,
\end{equation}
where $w(x)$ is a weight function on $(a,b)$ and $f$ is a smooth function. Let $r\in\mathbb{N}$, if the value of $f(x)$ and its consecutive derivatives up to an order $r-1$ are prescribed at the left endpoint, then an interpolatory quadrature rule can be constructed as
\begin{align}\label{eq:GRQ}
Q_{n,r}^{\mathrm{GR}}(f) = \sum_{j=0}^{r-1} w_j^{a} f^{(j)}(a) +
\sum_{j=1}^{n} w_j f(x_j).
\end{align}
If the nodes $\{x_j\}_{j=1}^{n}$ and the weights
$\{w_j^a\}_{j=0}^{r-1}\cup\{w_j\}_{j=1}^{n}$ in \eqref{eq:GRQ} are
chosen to maximize the degree of exactness, i.e., $(I-Q_{n,r}^{\mathrm{GR}})(f)=0$ for $f\in\mathcal{P}_{2n+r-1}$,
then the quadrature rule \eqref{eq:GRQ} is known as a generalized Gauss-Radau quadrature rule\footnote{Strictly speaking, the rule $Q_{n,r}^{\mathrm{GR}}(f)$ is referred to as a Gauss-Radau quadrature rule when $r=1$ and a generalized Gauss-Radau quadrature rule when $r\geq2$. Here, we refer to $Q_{n,r}^{\mathrm{GR}}(f)$ as a generalized Gauss-Radau quadrature rule for all $r\geq1$ for the sake of simplicity.}.

Generalized Gauss-Radau quadrature rule has received certain attention in the past two decades and it is now known that the nodes $\{x_j\}_{j=1}^{n}$ are located in $(a,b)$ and the weights $\{w_j^{a}\}_{j=0}^{r-1}\cup\{w_j\}_{j=1}^{n}$ are all positive (see, e.g., \cite{Gautschi2004a,Gautschi2009,Joulak2009}). In the following we briefly describe the implementation of the generalized Gauss-Radau quadrature rule. Let $\{\psi_k\}_{k=0}^{\infty}$ be the sequence of monic polynomials
that are orthogonal with respect to the new weight function $w_r(x)=w(x)(x-a)^r$ on
$(a,b)$ and
\begin{equation}
\int_{a}^{b} w_r(x) \psi_k(x) \psi_j(x) \mathrm{d}x =
\gamma_k\delta_{k,j}, \nonumber
\end{equation}
where $\delta_{k,j}$ is the Kronecker delta and $\gamma_k>0$. Moreover, let $\{x_j^\mathrm{G},w_j^\mathrm{G}\}_{j=1}^{n}$ denote the nodes and weights of an $n$-point Gaussian quadrature with respect to the weight $w_r(x)$ on $(a,b)$, i.e.,
\begin{equation}\label{eq:GQ}
\int_{a}^{b} w_r(x) f(x)\mathrm{d}x = \sum_{j=1}^{n}
w_j^\mathrm{G} f(x_j^\mathrm{G}), \quad \forall{f}\in\mathcal{P}_{2n-1}.
\end{equation}
From \cite[Theorem 3.9]{Gautschi2004b} we know that
\begin{equation}
x_j = x_j^\mathrm{G}, \quad w_j = \frac{w_j^\mathrm{G}}{(x_j-a)^r}, \quad
j=1,\ldots,n, \nonumber
\end{equation}
and thus the interior nodes and weights of $Q_{n,r}^{\mathrm{GR}}(f)$, i.e., $\{x_j\}_{j=1}^{n}\cup\{w_j\}_{j=1}^{n}$, can be calculated from the nodes and weights of the Gaussian quadrature rule \eqref{eq:GQ}. As for the boundary weights of $Q_{n,r}^{\mathrm{GR}}(f)$, i.e., $\{w_j^a\}_{j=0}^{r-1}$, setting $f(x)=(x-a)^{i-1}\psi_n^2(x)$ with $i=1,2,\ldots,r$ in
\eqref{eq:GRQ} gives
\begin{equation}
\left(
  \begin{array}{ccccc}
    a_{11} & a_{12}   & {\cdots}& a_{1r}\\
      & a_{22}  & {\cdots}& a_{2r}\\
      &    & {\ddots} & {\vdots}\\
      &    &          & a_{rr}\\
  \end{array}
\right) \left(
  \begin{array}{cccc}
    w_0^{a} \\
    w_1^{a} \\
    {\vdots} \\
    w_{r-1}^{a} \\
  \end{array}
\right) = \left(
  \begin{array}{ccccc}
    b_1   \\
    b_2   \\
    {\vdots} \\
    b_r   \\
  \end{array}
\right), \nonumber
\end{equation}
where
\begin{align}
a_{ij} = \left[ (x-a)^{i-1} \psi_n^2(x) \right]_{x=a}^{(j-1)}, \quad
b_i = \int_{a}^{b} \omega(x) (x-a)^{i-1} \psi_n^2(x) \ud x.
\nonumber
\end{align}
Therefore, the boundary weights can be derived by solving the above upper triangular system. 
The diagonal elements are given by $a_{ii}=(i-1)!\psi_n^2(a)$, where $\psi_n$ is the monic orthogonal polynomial with respect to the weight function $w_r(x)=w(x)(x-a)^r$ on $(a,b)$. In the case of Jacobi and Laguerre polynomials, for instance, direct calculation with explicit formulas available shows that these quantities are extremely small or extremely large as $n$ grows, which could result in an underflow or overflow problem. However, this issue can be circumvented by rescaling the monic polynomials $\{\psi_n\}$ appropriately \cite{Gautschi2009}.

\subsection{Regularization of improper Riemann integrals}\label{sec:Abel}
In this subsection we introduce the regularization of improper Riemann integrals. For the improper Riemann integral $\int_{0}^{\infty} f(x) \ud x$, its Abel limit is defined by \cite[Chapter~4]{Wong2001}
\begin{equation}\label{eq:Abel}
\lim_{s\rightarrow0^{+}} \int_{0}^{\infty} f(x) e^{-sx} \ud x.
\end{equation}
When the improper Riemann integral exists, e.g., $f$ is absolutely integrable on $(0,\infty)$, then the Abel limit is simply the improper integral itself. However, the importance of the Abel limit is that it may exist for certain divergent integrals that do not exist in the classical sense. Typical examples include the following two integrals
\begin{equation}\label{eq:FourierImp}
\int_{0}^{\infty} x^{\mu} e^{\mathrm{i}\omega x} \ud x = \frac{\Gamma(\mu+1)}{\omega^{\mu+1}} e^{(\mu+1)\pi\mathrm{i}/2},
\end{equation}
and
\begin{equation}\label{eq:BesselImp}
\int_{0}^{\infty} x^{\mu} J_{\nu}(\omega x) \ud x =
\frac{2^{\mu} \Gamma((\nu+\mu+1)/2)}{\omega^{\mu+1} \Gamma((\nu-\mu+1)/2)} ,
\end{equation}
where $\Gamma(z)$ is the gamma function and $\Re(\mu)>-1$ for \eqref{eq:FourierImp} and $\Re(\mu+\nu)>-1$ for \eqref{eq:BesselImp}. Note that \eqref{eq:FourierImp} and \eqref{eq:BesselImp} play an important role in studying the asymptotic behaviors of Fourier and Hankel transforms (see, e.g., \cite[Chapter 4]{Wong2001}). Moreover, \eqref{eq:BesselImp} was also used in \cite{Asheim2013} to compute the Gaussian quadrature rule with respect to the weight function $J_{\nu}(x)$ on $(0,\infty)$.

Below we state the first main result of this work, which gives conditions on the rotation of integration path of Fourier and Hankel transforms in the right-half plane. In the remainder of this paper, all improper integrals are defined using their Abel limits.
\begin{theorem}\label{thm:Rot}
If $f$ is analytic in the right half-plane and $|f(z)|\leq\mathcal{K}|z|^{\sigma}$ for some $\sigma\in\mathbb{R}$ as $z\rightarrow\infty$ and $|f(z)|\leq \mathcal{M}$ for some $\mathcal{M}>0$ as $z\rightarrow0$, then for $\Re(\mu)>-1$,
\begin{equation}\label{eq:FouRot}
\int_{0}^{\infty} f(x) x^{\mu} e^{\mathrm{i}\omega x} \ud x = e^{(\mu+1)\pi\mathrm{i}/2} \int_{0}^{\infty} f(\mathrm{i}x) x^{\mu} e^{-\omega x} \ud x,
\end{equation}
and for $\Re(\mu\pm\nu)>-1$,
\begin{equation}\label{eq:BesRot}
\int_{0}^{\infty} f(x) x^{\mu} J_{\nu}(\omega x) \ud x = \int_{0}^{\infty} \widehat{f}(\mathrm{i}x) x^{\mu} K_{\nu}(\omega x) \ud x,
\end{equation}
where $\widehat{f}(x)=(e^{(\mu-\nu)\pi\mathrm{i}/2}f(x)+e^{(\nu-\mu)\pi\mathrm{i}/2}f(-x))/\pi$ and $K_{\nu}(z)$ is the modified Bessel function of the second kind.
\end{theorem}
\begin{proof}
We only sketch the proof of \eqref{eq:BesRot} and the proof of \eqref{eq:FouRot} is similar. Let $C_{\varepsilon}$ with $\varepsilon>0$ denote the circle of radius $\varepsilon$ in the first quadrant and let $C_{\lambda}$ with $\lambda>\varepsilon$ denote the circle of radius $\lambda$ in the first quadrant. For $s>0$, by the identity \cite[Equation~(10.27.9)]{Frank2010handbook}, i.e., $\mathrm{i}\pi J_{\nu}(z) = e^{-\nu\pi\mathrm{i}/2} K_{\nu}(-\mathrm{i}z) - e^{\nu\pi\mathrm{i}/2} K_{\nu}(\mathrm{i}z)$ for $|\arg(z)|\leq\pi/2$, we have
\begin{align}\label{eq:P1}
\int_{\varepsilon}^{\lambda} e^{-sx} f(x) x^{\mu} J_{\nu}(\omega x) \ud x &= \frac{e^{-\nu\pi\mathrm{i}/2}}{\mathrm{i}\pi} \int_{\varepsilon}^{\lambda} e^{-sx} f(x) x^{\mu} K_{\nu}(-\mathrm{i}\omega x) \ud x \nonumber \\
&~~ - \frac{e^{\nu\pi\mathrm{i}/2}}{\mathrm{i}\pi} \int_{\varepsilon}^{\lambda} e^{-sx} f(x) x^{\mu} K_{\nu}(\mathrm{i}\omega x) \ud x.
\end{align}
For the first integral on the right-hand side, by Cauchy's theorem, we know that
\begin{align}\label{eq:P2}
\int_{\varepsilon}^{\lambda} e^{-sx} f(x) x^{\mu} K_{\nu}(-\mathrm{i}\omega x) \ud x &= \left(\int_{C_{\varepsilon}} - \int_{C_{\lambda}} \right) e^{-sz} f(z) z^{\mu} K_{\nu}(-\mathrm{i}\omega z) \ud z \nonumber \\
&~ + e^{(\mu+1)\pi\mathrm{i}/2} \int_{\varepsilon}^{\lambda} e^{-\mathrm{i}sx} f(\mathrm{i}x) x^{\mu} K_{\nu}(\omega x) \ud x,
\end{align}
where the contours are taken in the counterclockwise direction. From \cite[Chapter 10]{Frank2010handbook} we know that $K_{-\nu}(z)=K_{\nu}(z)$, and $K_{\nu}(z)=O(z^{-\nu})$ for $\Re(\nu)>0$ and $K_0(z)=O(\ln z)$ as $z\rightarrow0$. In the case of $\Re(\nu)=0$ and $\nu\neq0$, $K_{\nu}(z)=O(1)$ as $z\rightarrow0$. Hence,
\begin{align}
\left| \int_{C_{\varepsilon}} e^{-sz} f(z) z^{\mu} K_{\nu}(-\mathrm{i}\omega z) \ud z \right| &= \left\{\begin{array}{ll}
       O(\varepsilon^{\Re(\mu+1-\nu)}),          & \Re(\nu)>0, \\[1ex]
       O(\varepsilon^{\Re(\mu+1)}\ln\varepsilon),   & \Re(\nu)=0,~\Im(\nu)=0, \\[1ex]
       O(\varepsilon^{\Re(\mu+1)}),                 & \Re(\nu)=0,~\Im(\nu)\neq0, \\[1ex]
       O(\varepsilon^{\Re(\mu+1+\nu)}),          & \Re(\nu)<0,
              \end{array}
              \right. \nonumber
\end{align}
and therefore the contour integral on the left-hand side vanishes as $\varepsilon\rightarrow0^{+}$. Moreover, parametrizing $C_{\lambda}$ by $z=\lambda e^{\mathrm{i}\theta}$ with $0\leq\theta\leq\pi/2$ and using \cite[Equation~(10.25.3)]{Frank2010handbook}, we have
\begin{align}
\left| \int_{C_{\lambda}} e^{-sz} f(z) z^{\mu} K_{\nu}(-\mathrm{i}\omega z) \ud z \right| &\leq \mathcal{K} \lambda^{\mu+\sigma+1} \int_{0}^{\pi/2} e^{-s\lambda\cos\theta} |K_{\nu}(-\mathrm{i}\omega \lambda e^{\mathrm{i}\theta})| \ud \theta \nonumber \\
&\leq \mathcal{K} \lambda^{\mu+\sigma+1/2}
\int_{0}^{\pi/2} e^{-\lambda(s\cos\theta+\omega\sin\theta)} \ud \theta, \nonumber \\
&\leq \mathcal{K} \lambda^{\mu+\sigma+1/2} e^{-\lambda s}, \nonumber
\end{align}
where we have used the fact that the maximum of the integrand in the second inequality is attained at $\theta=0$ for $\omega\geq s>0$. Therefore, the contour integral on the left-hand side vanishes as $\lambda\rightarrow\infty$.
Letting $\varepsilon\rightarrow0^{+}$ and $\lambda\rightarrow\infty$ in \eqref{eq:P2}, we obtain
\begin{align}
\int_{0}^{\infty} e^{-sx} f(x) x^{\mu} K_{\nu}(-\mathrm{i}\omega x) \ud x &= e^{(\mu+1)\pi\mathrm{i}/2} \int_{0}^{\infty} e^{-\mathrm{i}sx} f(\mathrm{i}x) x^{\mu} K_{\nu}(\omega x) \ud x. \nonumber
\end{align}
and similarly,
\begin{align}
\int_{0}^{\infty} e^{-sx} f(x) x^{\mu} K_{\nu}(\mathrm{i}\omega x) \ud x &= e^{-(\mu+1)\pi\mathrm{i}/2} \int_{0}^{\infty} e^{\mathrm{i}sx} f(-\mathrm{i}x) x^{\mu} K_{\nu}(\omega x) \ud x. \nonumber
\end{align}
Combining the above two equations with \eqref{eq:P1} and letting $s\rightarrow0^{+}$ gives \eqref{eq:BesRot}. This ends the proof.
\end{proof}

We have the following remarks to Theorem \ref{thm:Rot}.
\begin{remark}
Both \eqref{eq:FourierImp} and \eqref{eq:BesselImp} can also be simply derived by setting $f(x)=1$ in Theorem \ref{thm:Rot}. Note that the condition in \eqref{eq:BesRot} is more restrictive than the condition in \eqref{eq:BesselImp}, which implies that the condition $\Re(\mu\pm\nu)>-1$ in \eqref{eq:BesRot} might be relaxed.
\end{remark}

\begin{remark}
Equation \eqref{eq:BesRot} was actually used in \cite[Theorem 3.5]{Asheim2013} to construct Gaussian quadrature rule with respect to $J_0(x)$ on $(0,\infty)$ with the requirement that $f$ is analytic in the right half-plane and suitably decaying at infinity. However, the validity of this requirement was not proved therein.
\end{remark}

\section{Complex generalized Gauss-Radau quadrature rules for Hankel transform of integer order}\label{sec:HankelInteger}
In this section, we consider the construction of complex generalized Gauss-Radau quadrature rules for the Hankel transform \eqref{eq:HankelTrans} of integer order, i.e., $\nu\in\mathbb{Z}$. Note that since $J_{-\nu}(z)=(-1)^{\nu}J_{\nu}(z)$, it is enough to consider the case $\nu\in\mathbb{N}_0$. We point out that the condition $\nu\in\mathbb{Z}$ might be relaxed when discussing the existence of orthogonal polynomials introduced in \eqref{eq:PolyBesGen} and we will specify it clearly.

For $\nu\in\mathbb{R}$, let $\mu+\nu>-1$ and $\mu-\nu\in\mathbb{N}_0$ and let
\begin{equation}\label{def:WeightFun}
w_{\mu,\nu}(x) = \frac{K_{\nu}(\sqrt{x})}{2} \left\{\begin{array}{ll}
       x^{(\mu-1)/2},  &  \mbox{$\mu-\nu$ even}, \\[1ex]
       x^{\mu/2},      &  \mbox{$\mu-\nu$ odd}.
              \end{array}
              \right.
\end{equation}
Recall that $K_{\nu}(\sqrt{x})>0$ for $x>0$ and $K_{\nu}(\sqrt{x}) \sim (\pi/2)^{1/2} x^{-1/4} e^{-\sqrt{x}}$ as $x\rightarrow\infty$ (see, e.g., \cite[Chapter~10]{Frank2010handbook}), and thus $w_{\mu,\nu}(x)$ is a weight function on $(0,\infty)$. Let $\{\phi_n\}_{n=0}^{\infty}$ denote the sequence of monic polynomials that are orthogonal with respect to the weight $w_{\mu,\nu}(x)$ and
\begin{equation}\label{eq:phi}
\int_{0}^{\infty} w_{\mu,\nu}(x) \phi_n(x) \phi_m(x) \ud x = \left\{\begin{array}{ll}
       0,         &  n\neq m, \\[1ex]
       \tau_n,  &  n=m.
              \end{array}
              \right.
\end{equation}
The $n$-point Gaussian quadrature rule with respect to $w_{\mu,\nu}(x)$ is
\begin{equation}\label{def:GaussK}
\int_{0}^{\infty} w_{\mu,\nu}(x) f(x) \ud x = \sum_{j=1}^{n} w_j f(x_j), \quad \forall f\in\mathcal{P}_{2n-1},
\end{equation}
where $\{x_j\}_{j=1}^{n}$ are the zeros of $\phi_n(x)$ and $w_j = \tau_{n-1}/(\phi_n{'}(x_j) \phi_{n-1}(x_j))$. By the properties of orthogonal polynomials and Gaussian quadrature rules, we know that $x_j\in(0,\infty)$ and $w_j>0$ for $j=1,\ldots,n$.

\begin{remark}
Up to a scaling factor, the weight function $w_{\mu,\nu}(x)$ is a special case of the generalized Prudnikov weight function $w_{\nu}^{\alpha}(x)=2x^{\alpha+\nu/2}K_{\nu}(2\sqrt{x})$, where $\nu\geq0$ and $\alpha>-1$, and thus the sequence of orthogonal polynomials $\{\phi_n\}_{n=0}^{\infty}$ is a special case of the generalized Prudnikov polynomials \cite{Gautschi2022,Yakubovich2021}.
\end{remark}

Below we state the second main result of this work.
\begin{theorem}\label{them:HankelInteger}
Suppose that $f$ is analytic in the right half-plane and
$|f(z)|\leq\mathcal{K}|z|^{\sigma}$ for some $\sigma\in\mathbb{R}$ as $z\rightarrow\infty$ and suppose that $f$ is analytic at $z=0$. Let $\{x_j,w_j\}_{j=1}^{n}$ be the nodes and weights of the Gaussian quadrature defined in \eqref{def:GaussK} and let $(\mathcal{Q}_{2n,\mu}^{\mathrm{HI}}f)(\omega)$ denote the quadrature rule of the form
\begin{equation}\label{eq:GGaussRadau}
(\mathcal{Q}_{2n,\mu}^{\mathrm{HI}}f)(\omega) = \frac{1}{\omega} \left( \sum_{k=0}^{\mu-1}
\frac{\hat{w}_{k}^{0}}{\omega^{k}} f^{(k)}(\hat{x}_0) + \sum_{j=1}^{2n} \hat{w}_j
f\left(\frac{\hat{x}_j}{\omega} \right) \right),
\end{equation}
where $\{\hat{x}_j\}_{j=1}^{2n} = \{\pm\mathrm{i}\sqrt{x_j} \}_{j=1}^{n}$ and $\hat{x}_0=0$, and
\begin{equation}
\{ \hat{w}_j \}_{j=1}^{2n} = \left\{\displaystyle \exp\left(\mp\frac{\nu\pi\mathrm{i}}{2}\right) \frac{w_j x_j^{-\kappa/2}}{\pi} \right\}_{j=1}^{n}, \nonumber
\end{equation}
and for $k=0,\ldots,\mu-1$,
\begin{equation}
\hat{w}_k^0 = \frac{1}{k!} \left( 2^{k} \frac{\Gamma((\nu+k+1)/2)}{\Gamma((\nu-k+1)/2)} - \frac{2}{\pi}\cos\left(\frac{(k-\nu)\pi}{2}\right)
\sum_{j=1}^{n} w_j x_j^{(k-\kappa)/2} \right), \nonumber
\end{equation}
and $\kappa=\mu$ when $\mu-\nu$ is even and $\kappa=\mu+1$ when $\mu-\nu$ is odd.
Then, for the Hankel transform \eqref{eq:HankelTrans} with $\nu\in\mathbb{N}_0$, the quadrature rule \eqref{eq:GGaussRadau} is exact for $f\in\mathcal{P}_{4n+\mu-1}$ when $\mu-\nu$ is even and for $f\in\mathcal{P}_{4n+\mu}$ when $\mu-\nu$ is odd.
\end{theorem}
\begin{proof}
Let $\mathcal{T}_{\mu}$ denote the Taylor expansion of $f$ of
degree $\mu-1$ at $\hat{x}_0 = 0$. It follows that
\begin{align}
(\mathcal{H}_{\nu}f)(\omega) & = \int_{0}^{\infty} \mathcal{T}_{\mu}(x)
J_{\nu}(\omega x) \ud x + \int_{0}^{\infty} \mathcal{R}_{\mu}(x) x^{\mu} J_{\nu}(\omega x) \ud x, \nonumber
\end{align}
where $\mathcal{R}_{\mu}(x) = (f(x)-\mathcal{T}_{\mu}(x))/x^{\mu}$. For the first integral on the right-hand side, using \eqref{eq:BesselImp} we have
\begin{align}
\int_{0}^{\infty} \mathcal{T}_{\mu}(x)
J_{\nu}(\omega x) \ud x &= \sum_{k=0}^{\mu-1} \frac{f^{(k)}(\hat{x}_0)}{k!} \int_{0}^{\infty} x^k J_{\nu}(\omega x) \ud x = \sum_{k=0}^{\mu-1} \frac{f^{(k)}(\hat{x}_0) 2^{k}}{k! \omega^{k+1}} \frac{\Gamma((\nu+k+1)/2)}{\Gamma((\nu-k+1)/2)}. \nonumber
\end{align}
For the second integral, using Theorem \ref{thm:Rot} we have
\begin{equation}
\int_{0}^{\infty} \mathcal{R}_{\mu}(x) x^{\mu} J_{\nu}(\omega x) \ud x = \int_{0}^{\infty} \widehat{\mathcal{R}}_{\mu}(\mathrm{i}x) x^\mu K_{\nu}(\omega x) \ud x, \nonumber
\end{equation}
where $\widehat{\mathcal{R}}_{\mu}(x) = (e^{(\mu-\nu)\mathrm{i}\pi/2} \mathcal{R}_{\mu}(x) + e^{(\nu-\mu)\mathrm{i}\pi/2} \mathcal{R}_{\mu}(-x) )/\pi$. It is easily verified that $\widehat{\mathcal{R}}_{\mu}(x)$ is an even function when $\mu-\nu$ is even and is an odd function when $\mu-\nu$ is odd. By the parity of $\widehat{\mathcal{R}}_{\mu}(x)$ and using the transformation $x\mapsto\sqrt{x}/\omega$ to the integral on the right-hand side yields
\begin{align}\label{eq:SecIntGauss}
\int_{0}^{\infty} \mathcal{R}_{\mu}(x) x^{\mu} J_{\nu}(\omega x) \ud x 
&= \frac{1}{\omega^{\mu+1}} \left\{\begin{array}{ll}
      {\displaystyle \int_{0}^{\infty} \widehat{\mathcal{R}}_{\mu}\left(\frac{\mathrm{i}\sqrt{x}}{\omega}\right) w_{\mu,\nu}(x) \ud x },  &  \mbox{$\mu-\nu$ even}, \\[3ex]
      {\displaystyle \int_{0}^{\infty} \frac{1}{\sqrt{x}} \widehat{\mathcal{R}}_{\mu}\left(\frac{\mathrm{i}\sqrt{x}}{\omega}\right) w_{\mu,\nu}(x) \ud x },  &  \mbox{$\mu-\nu$ odd},
              \end{array}
              \right. \nonumber \\
&\approx \frac{1}{\omega^{\mu+1}} \left\{\begin{array}{ll}
      {\displaystyle \sum_{j=1}^{n} w_j \widehat{\mathcal{R}}_{\mu}\left(\frac{\mathrm{i}\sqrt{x_j}}{\omega}\right) },  &  \mbox{$\mu-\nu$ even}, \\[3ex]
      {\displaystyle \sum_{j=1}^{n} \frac{w_j}{\sqrt{x_j}}  \widehat{\mathcal{R}}_{\mu}\left(\frac{\mathrm{i}\sqrt{x_j}}{\omega}\right) },  &  \mbox{$\mu-\nu$ odd},
              \end{array}
              \right.
\end{align}
where we have used the $n$-point Gausssian quadrature rule in \eqref{def:GaussK} to evaluate the integrals in the first line. Combining all the above results and after some simplification gives the quadrature rule \eqref{eq:GGaussRadau}. Moreover, the Gaussian quadrature rule in \eqref{eq:SecIntGauss} is exact for $\widehat{\mathcal{R}}_{\mu}\in\mathcal{P}_{4n-2}$ when $\mu-\nu$ is even and for $\widehat{\mathcal{R}}_{\mu}\in\mathcal{P}_{4n-1}$ when $\mu-\nu$ is odd.

Finally, we show the exactness of the quadrature rule \eqref{eq:GGaussRadau}. When $\mu-\nu$ is even, by the Taylor expansion of $f$, we get
\begin{equation}
\mathcal{R}_{\mu}(x) = \sum_{j=0}^{\infty} \frac{f^{(j+\mu)}(0)}{(j+\mu)!} x^j
~~~ \Rightarrow ~~~  \widehat{\mathcal{R}}_{\mu}(x) = \frac{2}{\pi} e^{(\mu-\nu)\pi\mathrm{i}/2} \sum_{j=0}^{\infty}
\frac{f^{(2j+\mu)}(0)}{(2j+\mu)!} x^{2j}.  \nonumber
\end{equation}
Recall that the quadrature rule in \eqref{eq:SecIntGauss} is exact for $\widehat{\mathcal{R}}_{\mu}\in\mathcal{P}_{4n-2}$, we therefore deduce that the quadrature rule \eqref{eq:GGaussRadau} is exact for $f\in\mathcal{P}_{4n+\mu-1}$. When $\mu-\nu$ is odd, it can be shown in a similar way that the quadrature rule \eqref{eq:GGaussRadau} is exact for $f\in\mathcal{P}_{4n+\mu}$. This ends the proof.
\end{proof}

Some remarks on Theorem \ref{them:HankelInteger} are in order.
\begin{itemize}
\item[(i)] When $\nu$ is an even integer, the weights $\{ \hat{w}_j \}_{j=1}^{2n}$ are all real and the weights corresponding to the nodes $\{\mathrm{i}\sqrt{x_j}\}_{j=1}^{n}$ are the same as the weights corresponding to $\{-\mathrm{i}\sqrt{x_j}\}_{j=1}^{n}$. When $\nu$ is an odd integer, then the weights $\{ \hat{w}_j \}_{j=1}^{2n}$ are all purely imaginary and the weights corresponding to $\{\mathrm{i}\sqrt{x_j}\}_{j=1}^{n}$ are the negative values of the weights corresponding to $\{-\mathrm{i}\sqrt{x_j}\}_{j=1}^{n}$.

\item[(ii)] When $\mu=\nu=0$, the nodes and weights of the quadrature rule $(\mathcal{Q}_{2n,\mu}^{\mathrm{HI}}f)(\omega)$ are $\{\hat{x}_j\}_{j=1}^{2n} = \{\pm\mathrm{i}\sqrt{x_j}\}_{j=1}^{n}$ and $\{\hat{w}_j\}_{j=1}^{2n} = \{ w_j/\pi\}_{j=1}^{n}\cup\{ w_j/\pi \}_{j=1}^{n}$, and $\{x_j,w_j\}_{j=1}^{n}$ are the nodes and weights of the $n$-point Gaussian quadrature rule with respect to $w_{0,0}(x)=x^{-1/2}K_{0}(\sqrt{x})/2$ on $(0,\infty)$. In this case, the rule $(\mathcal{Q}_{2n,\mu}^{\mathrm{HI}}f)(\omega)$ is exactly the Gaussian quadrature rule derived in \cite[Theorem 3.5]{Asheim2013} (Note that $w_j/\pi$ in $\{\hat{w}_j\}_{j=1}^{2n}$ was mistakenly written as $w_j/2$ in \cite[Theorem~3.5]{Asheim2013}).

\item[(iii)] If $f(x)$ is even and $\nu$ is an odd integer or if $f(x)$ is odd and $\nu$ is an even integer, then from item (i) above we can deduce that the second sum on the right hand side of \eqref{eq:GGaussRadau} vanishes and thus the quadrature rule $(\mathcal{Q}_{2n,\mu}^{\mathrm{HI}}f)(\omega)$ can be simplified.

\item[(iv)] The rules $(\mathcal{Q}_{2n,\mu}^{\mathrm{HI}}f)(\omega)$ are generalized Gauss-Radau quadrature rules in the sense that their exactness for polynomials is maximized.

\item[(v)] We take $\mu\in\mathbb{N}_0$ since we have used the Taylor expansion of $f(x)$ at $x=0$ in the construction of the quadrature rule and $\mu$ plays the role of order of derivatives of $f(x)$ taken at $x=0$. If taking $\mu\in\mathbb{R}$ and using the fractional Taylor expansion of $f(x)$, then $\mathcal{R}_{\mu}(x)$ will involve a linear combination of both smooth and singular terms, which make it difficult to construct a Gaussian quadrature rule to evaluate the resulting integral.
\end{itemize}

Asymptotic error estimates of the quadrature rule \eqref{eq:GGaussRadau}
are given below.
\begin{theorem}\label{them:GGRQ}
Under the assumptions of Theorem \ref{them:HankelInteger}, we have
\begin{equation}\label{eq:AsymError}
(\mathcal{H}_{\nu} f)(\omega) - (\mathcal{Q}_{2n,\mu}^{\mathrm{HI}}f)(\omega) = \left\{\begin{array}{ll}
       \mathcal{O}(\omega^{-4n-\mu-1}),    &  \mbox{$\mu-\nu$ even}, \\[2ex]
       \mathcal{O}(\omega^{-4n-\mu-2}),    &  \mbox{$\mu-\nu$ odd},
              \end{array}
              \right. \quad \omega\rightarrow\infty.
\end{equation}
\end{theorem}
\begin{proof}
Note that the rule $(\mathcal{Q}_{2n,\mu}^{\mathrm{HI}}f)(\omega)$ is exact for $f\in\mathcal{P}_{4n+\mu-1}$ when $\mu-\nu$ is even and for $f\in\mathcal{P}_{4n+\mu}$ when $\mu-\nu$ is odd, the asymptotic error estimates \eqref{eq:AsymError} follows immediately from \cite[Lemma 1.6]{Asheim2013}.
\end{proof}

\begin{remark}
For the rules $(\mathcal{Q}_{2n,\mu}^{\mathrm{HI}}f)(\omega)$ with $\mu=\nu+2k-1$ and $\mu=\nu+2k$ and $k\geq1$, it is easily verified that their nodes $\{\hat{x}_j\}_{j=1}^{2n}$ and weights $\{\hat{w}_j\}_{j=1}^{2n}$ are the same. Moreover, direct calculation shows that $\hat{w}_{\nu+2j-1}^0=0$ for $j\geq1$, and therefore the boundary weights $\{\hat{w}_j^0\}_{j=0}^{\nu+2k-2}$ for these two rules are also the same. We conclude that the rules $(\mathcal{Q}_{2n,\mu}^{\mathrm{HI}}f)(\omega)$ with $\mu-\nu=2k-1$ and $\mu-\nu=2k$ and $k\geq1$ are always the same and this explains why their asymptotic error estimates in \eqref{eq:AsymError} are the same.
\end{remark}

From Theorem \ref{them:HankelInteger} we see that the nodes of the generalized Gauss-Radau quadrature rule can be derived from the nodes of the Gaussian quadrature rule in \eqref{def:GaussK}. This implies the existence of the polynomials defined in \eqref{eq:PolyBesGen} for even degrees. Below we state a more general result, which gives the polynomials of all degrees when $\mu-\nu$ is even. When $\mu-\nu$ is odd, however, we show that the polynomials exist only for even degrees.
\begin{theorem}\label{thm:PolyConnection}
Let $\mu+\nu>-1$ and $\mu-\nu\in\mathbb{N}_0$ with $\nu\in\mathbb{R}$ and let $\{\phi_n\}_{n=0}^{\infty}$ be the polynomials defined in \eqref{eq:phi}. When $\mu-\nu$ is even, then for each $n\geq0$,
\begin{equation}\label{eq:EvenDegree}
P_{2n}^{(\mu,\nu)}(x) = (-1)^n \phi_n(-x^2),
\end{equation}
and
\begin{equation}\label{eq:OddDegree}
P_{2n+1}^{(\mu,\nu)}(x) = \frac{(-1)^{n+1}}{x} \left( \phi_{n+1}(-x^2) - \frac{\phi_{n+1}(0)}{\phi_n(0)} \phi_n(-x^2) \right).
\end{equation}
When $\mu-\nu$ is odd, then \eqref{eq:EvenDegree} still holds, but $P_{2n+1}^{(\mu,\nu)}(x)$ does not exist.
\end{theorem}
\begin{proof}
We first consider the case when $\mu-\nu$ is even. We only consider the proof of \eqref{eq:EvenDegree}, since the proof of \eqref{eq:OddDegree} is similar.
For any $s\in\mathcal{P}_{2n-1}$, by the definition of $P_{2n}^{(\mu,\nu)}(x)$,
\begin{equation}
\int_{0}^{\infty} P_{2n}^{(\mu,\nu)}(x) s(x) x^{\mu} J_{\nu}(x) \ud x = 0. \nonumber
\end{equation}
For the integral on the left-hand side, let $\Psi(x)=P_{2n}^{(\mu,\nu)}(x) s(x)$ and using Theorem \ref{thm:Rot} we have
\begin{align}
\int_{0}^{\infty} \Psi(x) x^{\mu} J_{\nu}(x) \ud x &= \frac{e^{(\mu-\nu)\mathrm{i}\pi/2}}{\pi} \int_{0}^{\infty} \left[ \Psi(\mathrm{i}x) + \Psi(-\mathrm{i}x) \right] x^{\mu} K_{\nu}(x) \ud x \nonumber \\
&= \frac{e^{(\mu-\nu)\mathrm{i}\pi/2}}{\pi} \int_{0}^{\infty} \left[ \Psi(\mathrm{i}\sqrt{x}) + \Psi(-\mathrm{i}\sqrt{x}) \right] w_{\mu,\nu}(x) \ud x, \nonumber
\end{align}
where we have used the transformation $x\mapsto\sqrt{x}$ in the last equality. If we set $s(x)=x^{2k+1}$, where $k=0,\ldots,n-1$, then
\begin{align}
\int_{0}^{\infty} \left[\sqrt{x} \left( P_{2n}^{(\mu,\nu)}(\mathrm{i}\sqrt{x}) - P_{2n}^{(\mu,\nu)}(-\mathrm{i}\sqrt{x}) \right) \right] x^{k} w_{\mu,\nu}(x) \ud x &= 0. \nonumber
\end{align}
Note that the term inside the square brackets is a polynomial of degree $n$ and is orthogonal to all polynomial of lower degree with respect to $w_{\mu,\nu}(x)$, we can deduce that
\[
\sqrt{x} \left( P_{2n}^{(\mu,\nu)}(\mathrm{i}\sqrt{x}) - P_{2n}^{(\mu,\nu)}(-\mathrm{i}\sqrt{x}) \right) = \lambda \phi_n(x),
\]
where $\lambda$ is a constant. Setting $x=0$ and noting that $\phi_n(0)\neq0$, we obtain $\lambda=0$, and thus $P_{2n}^{(\mu,\nu)}(x)$ is even. If we set $s(x)=x^{2k}$, where $k=0,\ldots,n-1$, then
\begin{align}
\int_{0}^{\infty} \left[ P_{2n}^{(\mu,\nu)}(\mathrm{i}\sqrt{x}) + P_{2n}^{(\mu,\nu)}(-\mathrm{i}\sqrt{x}) \right] x^k w_{\mu,\nu}(x) \ud x &= 0. \nonumber
\end{align}
Note that the term inside the square brackets is a polynomial of degree $n$ and is orthogonal to all polynomial of lower degree with respect to $w_{\mu,\nu}(x)$, we can deduce that
\begin{align}
P_{2n}^{(\mu,\nu)}(\mathrm{i}\sqrt{x}) + P_{2n}^{(\mu,\nu)}(-\mathrm{i}\sqrt{x}) = 2 (-1)^n \phi_n(x). \nonumber
\end{align}
Recall that $P_{2n}^{(\mu,\nu)}(x)$ is even, it follows that $P_{2n}^{(\mu,\nu)}(x) = (-1)^n \phi_n(-x^2)$. This proves the case where $\mu-\nu$ is even.

When $\mu-\nu$ is odd, \eqref{eq:EvenDegree} follows by similar arguments as above. Now we show that $P_{2n+1}^{(\mu,\nu)}(x)$ do not exist. By similar arguments as above we find that
\begin{equation}
\int_{0}^{\infty} \left[ \frac{P_{2n+1}^{(\mu,\nu)}(\mathrm{i}\sqrt{x}) - P_{2n+1}^{(\mu,\nu)}(-\mathrm{i}\sqrt{x})}{\sqrt{x}} \right] x^{k} w_{\mu,\nu}(x) \ud x = 0, \nonumber
\end{equation}
where $k=0,\ldots,n$. Note that the term inside the square brackets is a polynomial of degree $n$, we deduce that $P_{2n+1}^{(\mu,\nu)}(x)$ is an even function. However, this is impossible since its leading term is $x^{2n+1}$. We conclude that $P_{2n+1}^{(\mu,\nu)}(x)$ does not exist and this ends the proof.
\end{proof}

The following corollary can be derived from Theorem \ref{thm:PolyConnection} directly.
\begin{corollary}
When $\mu-\nu$ is even, then $P_n^{(\mu,\nu)}(x)$ always exists and is an even function when $n$ is even and is an odd function when $n$ is odd. When $\mu-\nu$ is odd, then $P_n^{(\mu,\nu)}(x)$ exists only for even $n$ and is an even function in this case. Moreover, in these cases that $P_n^{(\mu,\nu)}(x)$ exists, its zeros are all located on the imaginary axis and symmetric with respect to the real axis.
\end{corollary}

\begin{remark}
It is known that orthogonal polynomials can also be expressed in terms of the associated Hankel determinant \cite[Chapter 2]{Gautschi2004b}. In fact, when $\mu-\nu$ is odd, it is easily checked that the Hankel determinant associated with $P_n^{(\mu,\nu)}(x)$ vanishes for odd $n$, which also confirms the nonexistence of the polynomials with odd degrees.
\end{remark}

\begin{remark}
With the results in Theorem \ref{thm:PolyConnection}, we point out that the construction of complex generalized Gauss-Radau quadrature rule in Theorem \ref{them:HankelInteger} can also be extended to Hankel transform of the form
\[
\int_{0}^{\infty} x^{\alpha} f(x) J_{\nu}(\omega x) \mathrm{d}x,
\]
where $\alpha+\nu>-1$, $\alpha-\nu\in\mathbb{Z}$ and $f$ satisfies the same assumptions as in Theorem \ref{them:HankelInteger}.
\end{remark}

When implementing $(\mathcal{Q}_{2n,\mu}^{\mathrm{HI}}f)(\omega)$, one needs to calculate the Gaussian quadrature rule defined in \eqref{def:GaussK}. By \cite[Equation~(10.43.19)]{Frank2010handbook} we know that the moments of $w_{\mu,\nu}(x)$ can be written explicitly as
\begin{equation}
\int_{0}^{\infty} x^{k} w_{\mu,\nu}(x) \ud x = \left\{\begin{array}{ll}
       {\displaystyle \Gamma\left(k+\frac{\mu-\nu+1}{2}\right)
\Gamma\left(k+\frac{\mu+\nu+1}{2}\right) 2^{2k+\mu-1}  },    &  \mbox{$\mu-\nu$ even}, \\[3ex]
       {\displaystyle \Gamma\left(k+\frac{\mu-\nu+2}{2}\right)
\Gamma\left(k+\frac{\mu+\nu+2}{2}\right) 2^{2k+\mu}  },    &  \mbox{$\mu-\nu$ odd},
              \end{array}
              \right. \nonumber
\end{equation}
and the monic polynomials $\{\phi_n\}_{n=1}^{\infty}$ can be calculated by the Gram--Schmidt orthogonalization procedure. For example, explicit expressions of $\phi_1$ and $\phi_2$ are given below:
\begin{equation}\label{eq:PhiExpression1}
\phi_1(x) = \left\{\begin{array}{ll}
       {\displaystyle x - (\mu+1)^2+\nu^2},    &  \mbox{$\mu-\nu$ even}, \\[1ex]
       {\displaystyle x - (\mu+2)^2+\nu^2},    &  \mbox{$\mu-\nu$ odd},
              \end{array}
              \right.
\end{equation}
and
\begin{equation}\label{eq:PhiExpression2}
\phi_2(x) = \left\{\begin{array}{ll}
       {\displaystyle x^2 - 2b_{\mu,\nu}x + c_{\mu,\nu}},       &  \mbox{$\mu-\nu$ even}, \\[1ex]
       {\displaystyle x^2 - 2b_{\mu+1,\nu}x + c_{\mu+1,\nu} },  &  \mbox{$\mu-\nu$ odd},
              \end{array}
              \right.
\end{equation}
where $b_{\mu,\nu}$ and $c_{\mu,\nu}$ are given by
\begin{align}
b_{\mu,\nu} &= \frac{(\mu+3)(\mu-\nu+3)(\mu+\nu+3)}{\mu+2}, \nonumber \\
c_{\mu,\nu} &= \frac{(\mu+4)(\mu-\nu+3)(\mu+\nu+3)(\mu-\nu+1)(\mu+\nu+1)}{\mu+2}. \nonumber
\end{align}
Consequently, the nodes and weights of the Gaussian quadrature rule \eqref{def:GaussK} for $n=1,2$ can be calculated immediately using the above expressions. For larger $n$, the Gaussian quadrature rule \eqref{def:GaussK} can be calculated from the moments by the classical Chebyshev algorithm and high-precision arithmetic is required to reduce the underlying ill-conditioning (see \cite{Gautschi2004b,Gautschi2022}).

\section{Numerical examples}
In this section we present numerical experiments to confirm our findings and illustrate the performance of the proposed quadrature rule in \eqref{eq:GGaussRadau}. In our computations, the nodes and weights of the Gaussian quadrature rule \eqref{def:GaussK} are calculated by using \eqref{eq:PhiExpression1} and \eqref{eq:PhiExpression2} directly in Matlab with double precision for $n=1,2$ and by the moments of $w_{\mu,\nu}(x)$ in Maple with 100-digit arithmetic for larger $n$.

\subsection{Distribution of the zeros of $P_n^{(\mu,\nu)}(x)$}
\begin{figure}
\centering
\includegraphics[width=0.45\textwidth]{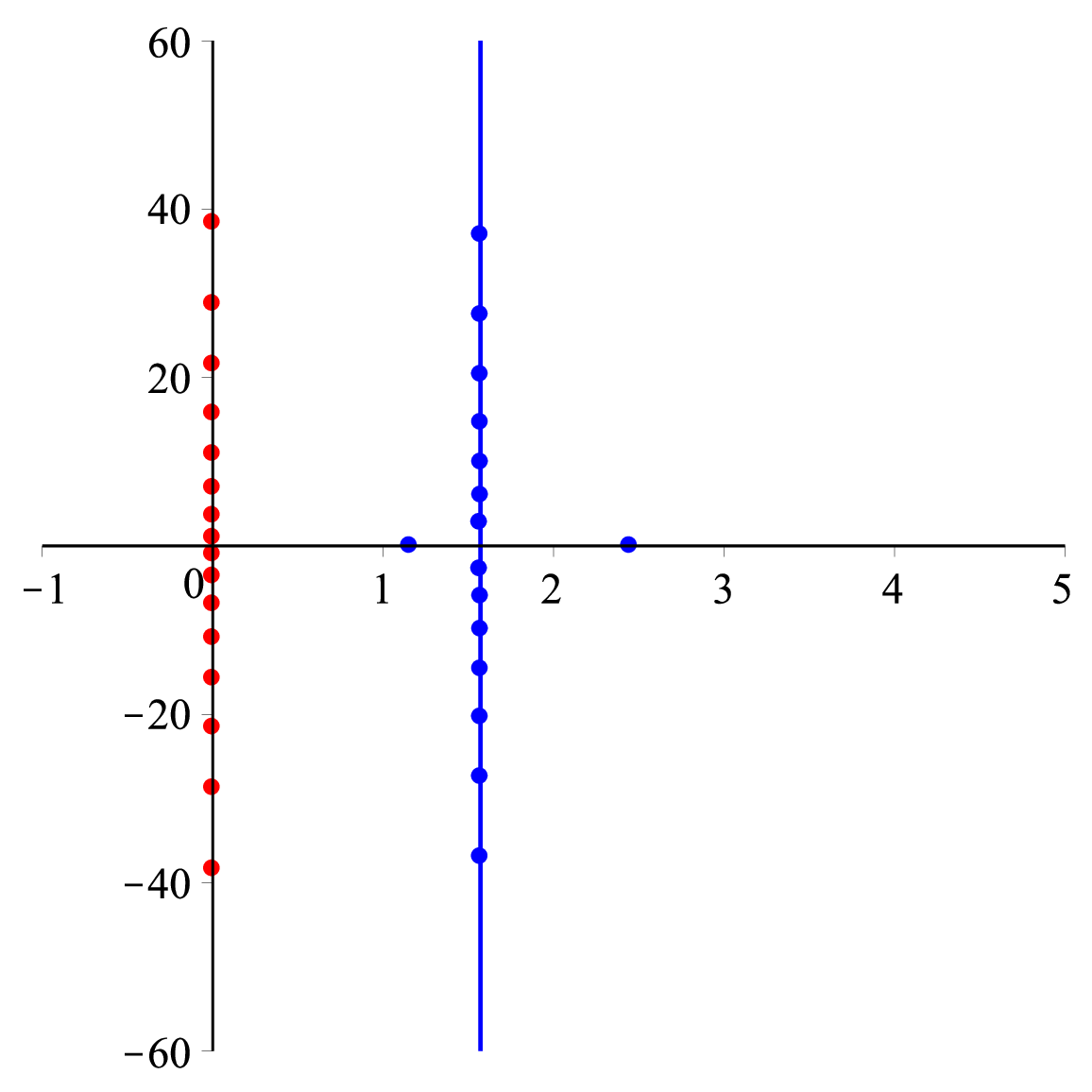}~~
\includegraphics[width=0.45\textwidth]{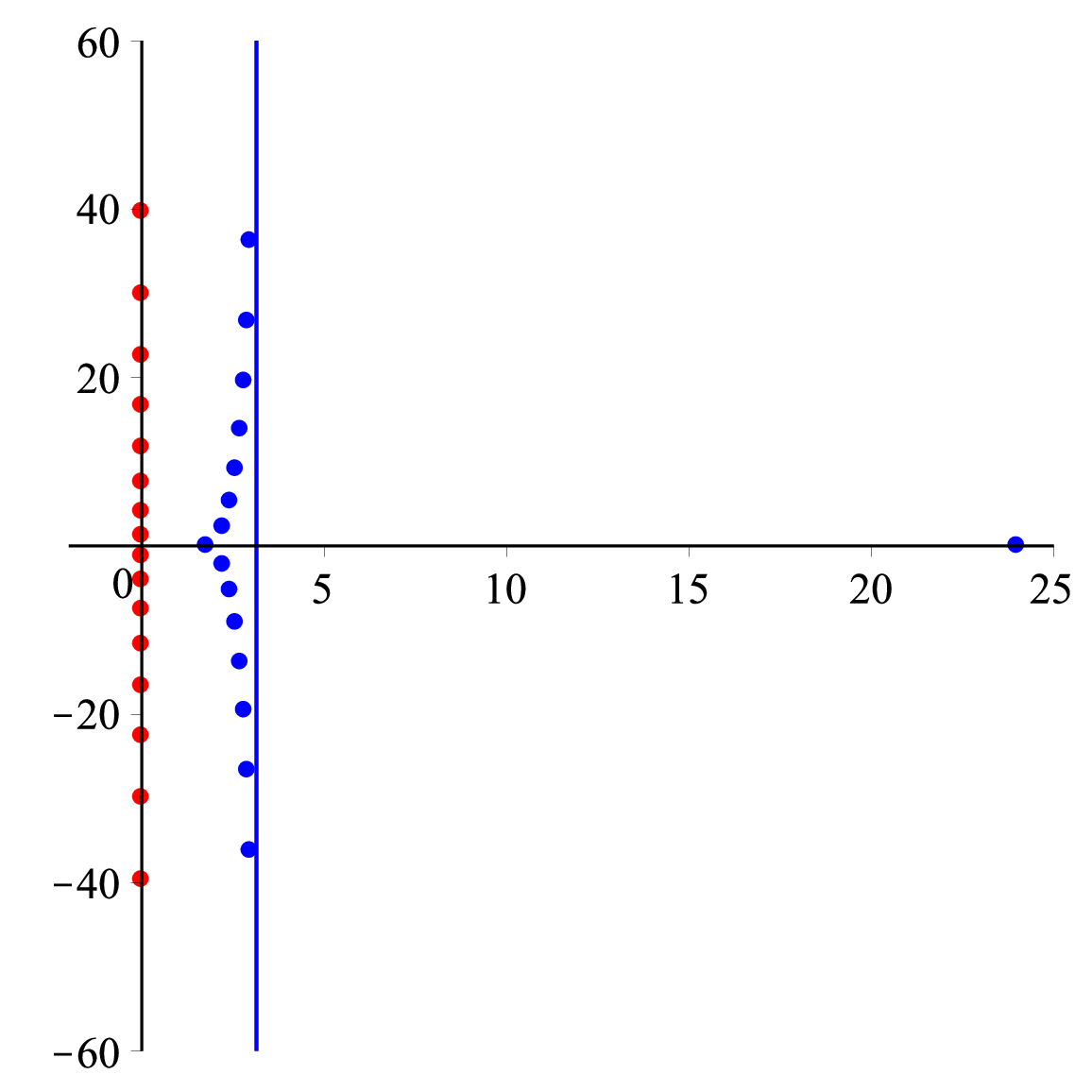}\\
\includegraphics[width=0.45\textwidth]{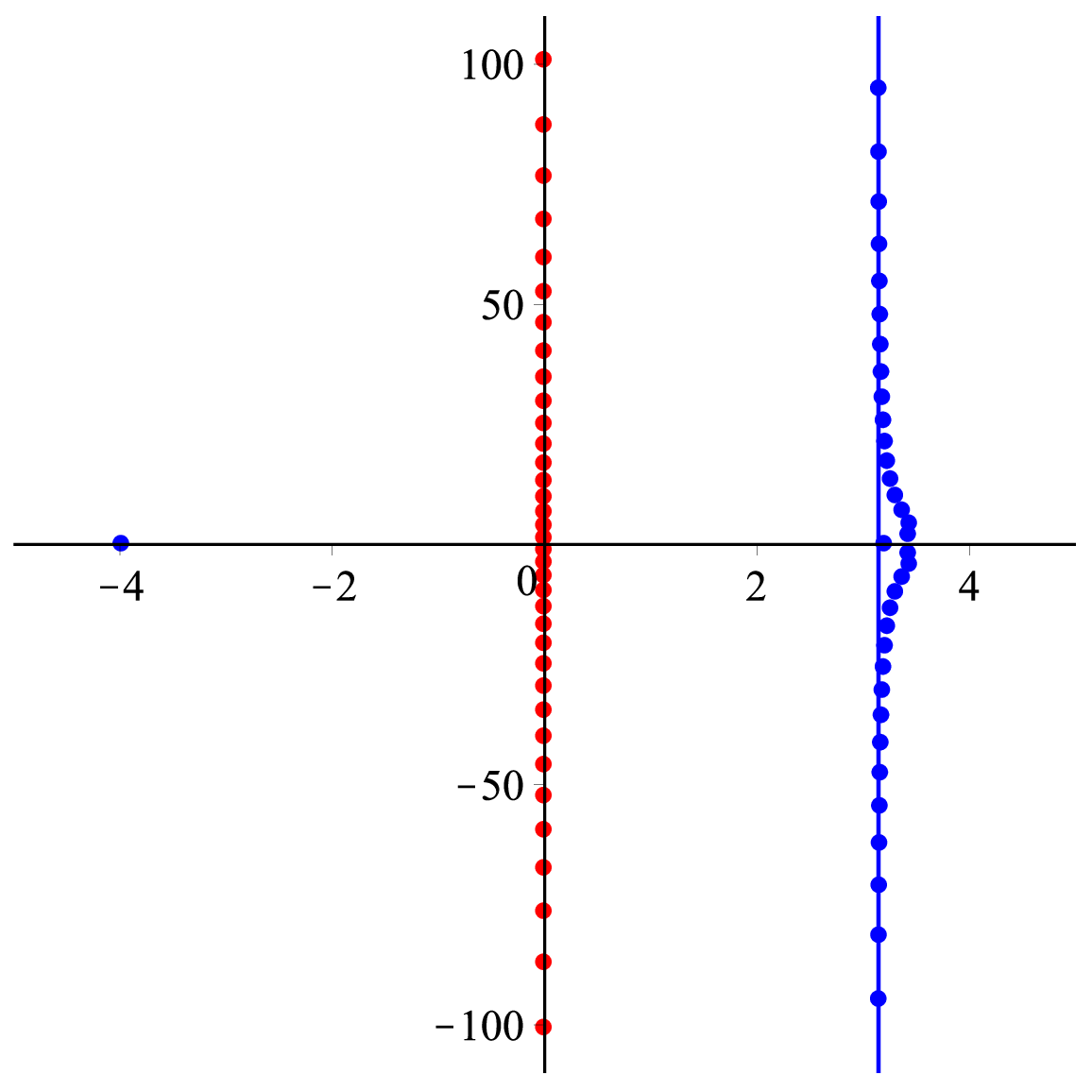}~~
\includegraphics[width=0.45\textwidth]{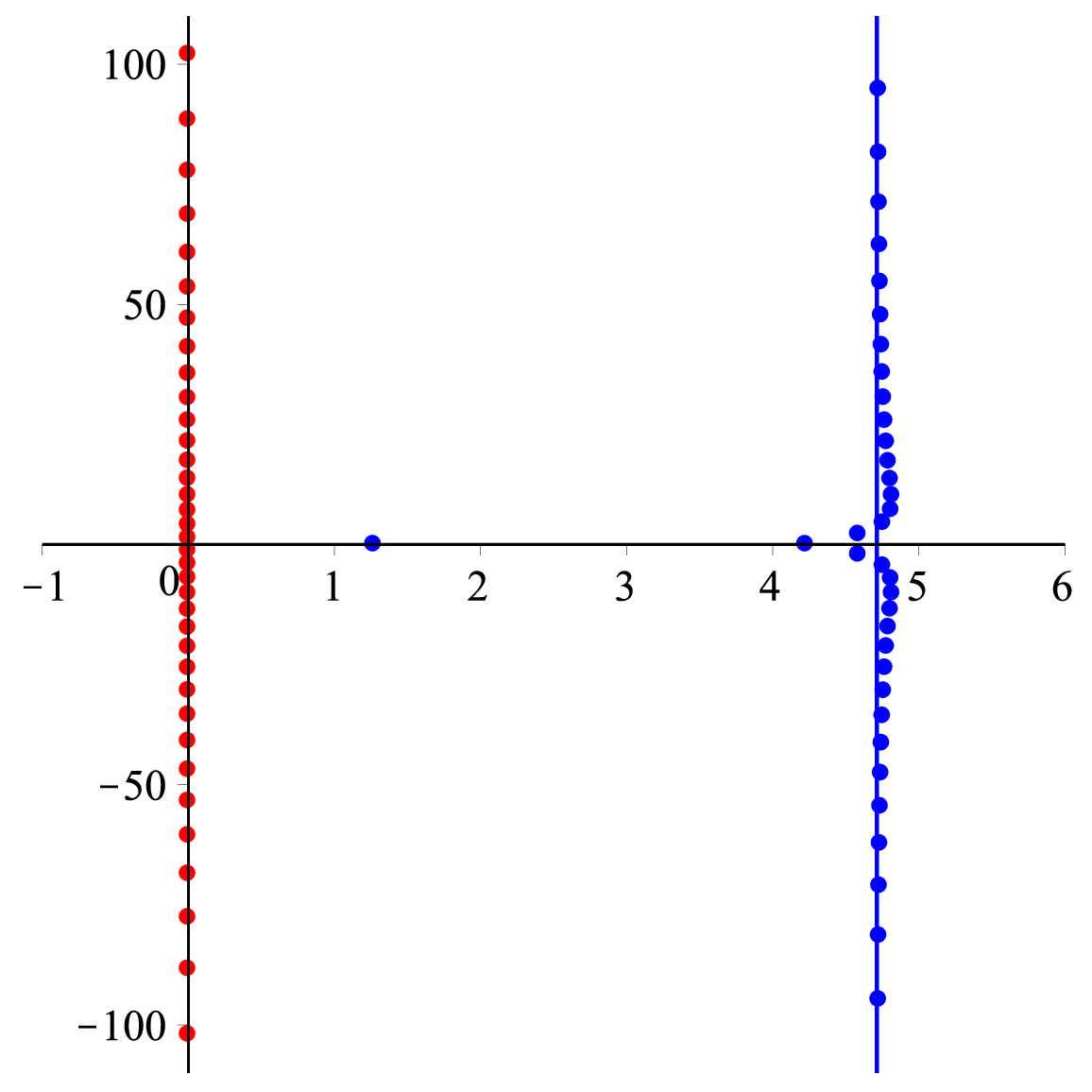}
\caption{The zeros of $P_{n}^{(0,\nu)}(x)$ (blue) and the zeros of
$P_{n}^{(\nu,\nu)}(x)$ (red). Top row shows $n=16$ for $\nu=1$ (left) and $\nu=2$ (right) and bottom row shows $n=36$ for $\nu=4$ (left) and $\nu=5$ (right). The vertical lines in the top row are $\Re(z)=\nu\pi/2$ and in the bottom row are $\Re(z)=(\nu-2)\pi/2$.} \label{fig:GRNodesInt}
\end{figure}

We consider the distribution of the zeros of $P_n^{(\mu,\nu)}(x)$.
In Figure \ref{fig:GRNodesInt} we display the zeros of $P_n^{(\mu,\nu)}(x)$ for four values of $\nu\in\mathbb{N}$ and we consider two choices of $\mu$: $\mu=0$ and $\mu=\nu$, which correspond respectively to the nodes of the Gaussian quadrature rule in \cite{Asheim2013} and the generalized Gauss-Radau quadrature rule in Theorem \ref{them:HankelInteger}. For $\mu=0$, Figure \ref{fig:GRNodesInt} shows that the zeros tend to cluster along the vertical line $\Re(z)=\nu\pi/2$ for $\nu=1,2$, but $\Re(z)=(\nu-2)\pi/2$ for $\nu=4,5$. In fact, the vertical line $\Re(z)=\nu\pi/2$ was already observed in \cite{Asheim2013} based on numerical calculations for $\nu=1/2,1,3/2,2$ and was further proved in \cite{Deano2016} for $\nu\in[0,1/2]$. However, we observe that there are transitions in the asymptotic distribution of the zeros and the vertical line becomes $\Re(z)=(\nu-2)\pi/2$ for $\nu\in[3,7]$ and will change further for larger $\nu$. On the other hand, we observe from Figure \ref{fig:GRNodesInt} that the nodes of the generalized Gauss-Radau quadrature rule are always located on the imaginary axis, as expected from Theorems \ref{them:HankelInteger} and \ref{thm:PolyConnection}.

\begin{figure}
\centering
\includegraphics[width=0.45\textwidth]{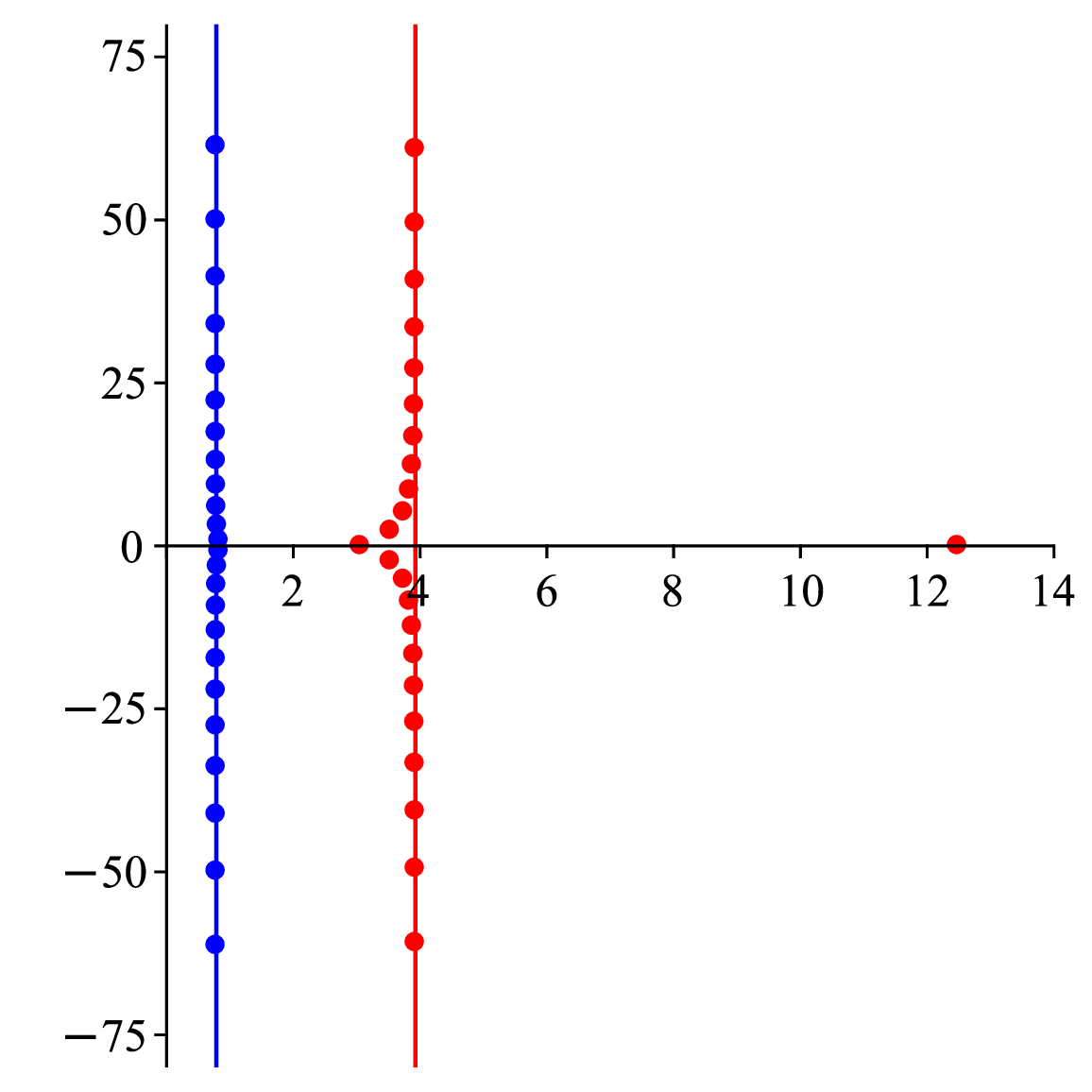}~~
\includegraphics[width=0.45\textwidth]{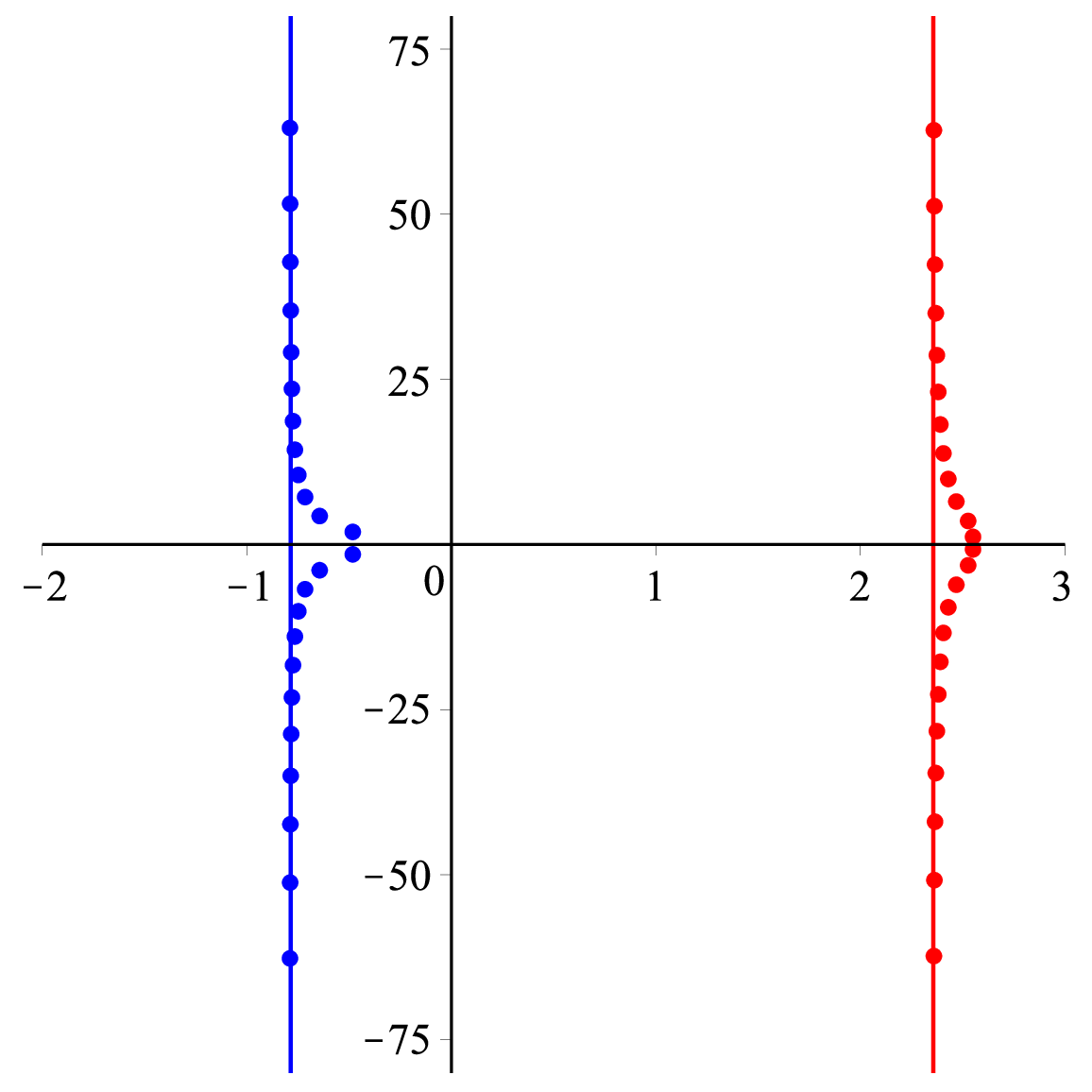}
\caption{The zeros of $P_n^{(\mu,\nu)}(x)$ for $\mu=1$ (left) and $\mu=2$ (right). Here $n=24$, $\nu=3/2$ (blue) and $\nu=7/2$ (red) and the vertical lines are $\Re(z)=(\nu-\mu)\pi/2$.} \label{fig:GRNodesFrac}
\end{figure}

What is the distribution of the zeros of $P_n^{(\mu,\nu)}(x)$ for $\nu\notin\mathbb{N}_0$? Figure \ref{fig:GRNodesFrac} shows the distribution of the zeros of $P_n^{(\mu,\nu)}(x)$ for $\nu=3/2,7/2$ and $\mu=1,2$. We see that the zeros of $P_n^{(\mu,\nu)}(x)$ tend to cluster along the vertical line $\Re(z)=(\nu-\mu)\pi/2$, but not the imaginary axis. Finally, to provide a more comprehensive insight into the distribution of the zeros of $P_n^{(\mu,\nu)}(x)$, we denote by $\Re(z)=\upsilon_{\mu,\nu}\pi/2$ the vertical line that the zeros of $P_n^{(\mu,\nu)}(x)$ cluster for large $n$ and we list in Table \ref{tab:line} the values of $\upsilon_{\mu,\nu}$ for several different values of $\mu$ and $\nu$. We see that $\upsilon_{\mu,\nu}$ changes both for the case of a fixed $\mu$ and increasing $\nu$ and for the case of a fixed $\nu$ and increasing $\mu$.

\begin{table}[h!]
  \begin{center}
    \caption{The value of $\upsilon_{\mu,\nu}$ for different values of $\mu$ and $\nu$.}
    \label{tab:line}
    \renewcommand\arraystretch{1.8}
    \begin{tabular}{c|cccccccccc}
      \hline
      \multicolumn{1}{c|}{\diagbox{$\mu$}{$\nu$}} &
      \multicolumn{1}{c}{$3/2$} &
      \multicolumn{1}{c}{$7/2$} &
      \multicolumn{1}{c}{$5$} &
      \multicolumn{1}{c}{$11/2$} &
      \multicolumn{1}{c}{$8$} &
      \multicolumn{1}{c}{$23/2$} \\
      \hline
{$1$} & {$\nu-\mu$} & {$\nu-\mu$} & {$\nu-\mu-2$} & {$\nu-\mu-2$} & {$\nu-\mu-2$} & {$\nu-\mu-4$} \\
{$2$} & {$\nu-\mu$} & {$\nu-\mu$} & {$\nu-\mu$} & {$\nu-\mu$} & {$\nu-\mu-2$} & {$\nu-\mu-4$} \\
{$3$} & {$\nu-\mu+2$} & {$\nu-\mu$}&{$\nu-\mu$}   & {$\nu-\mu$}   & {$\nu-\mu-2$}   & {$\nu-\mu-2$}\\
{$4$} & {$\nu-\mu+2$} & {$\nu-\mu$}&{$\nu-\mu$}   & {$\nu-\mu$}   & {$\nu-\mu$}   & {$\nu-\mu-2$}\\
{$5$} & {$\nu-\mu+4$} & {$\nu-\mu+2$}&{$\nu-\mu$} & {$\nu-\mu$}   & {$\nu-\mu$}   & {$\nu-\mu-2$}\\
      \hline
    \end{tabular}
  \end{center}
\end{table}

\subsection{Performance of the rule $(\mathcal{Q}_{2n,\mu}^{\mathrm{HI}}f)(\omega)$}
\begin{figure}
\centering
\includegraphics[width=7.4cm,height=6.2cm]{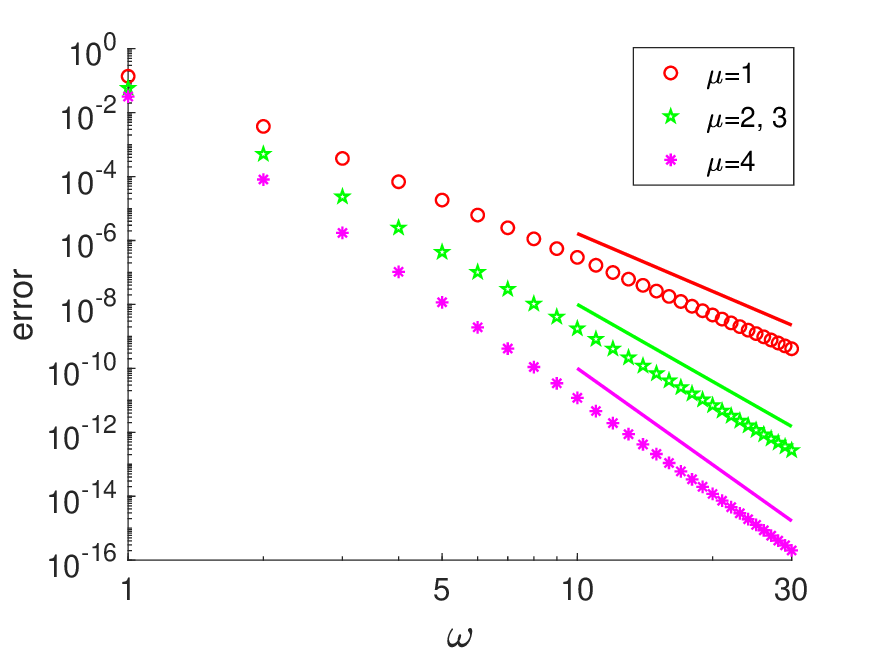}~~
\includegraphics[width=7.4cm,height=6.2cm]{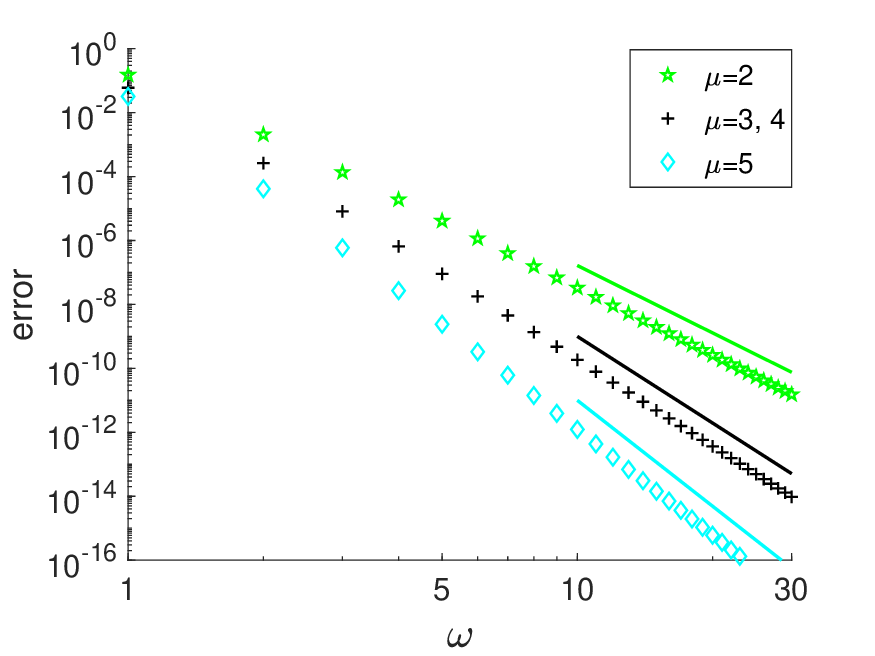}\\
\includegraphics[width=7.4cm,height=6.2cm]{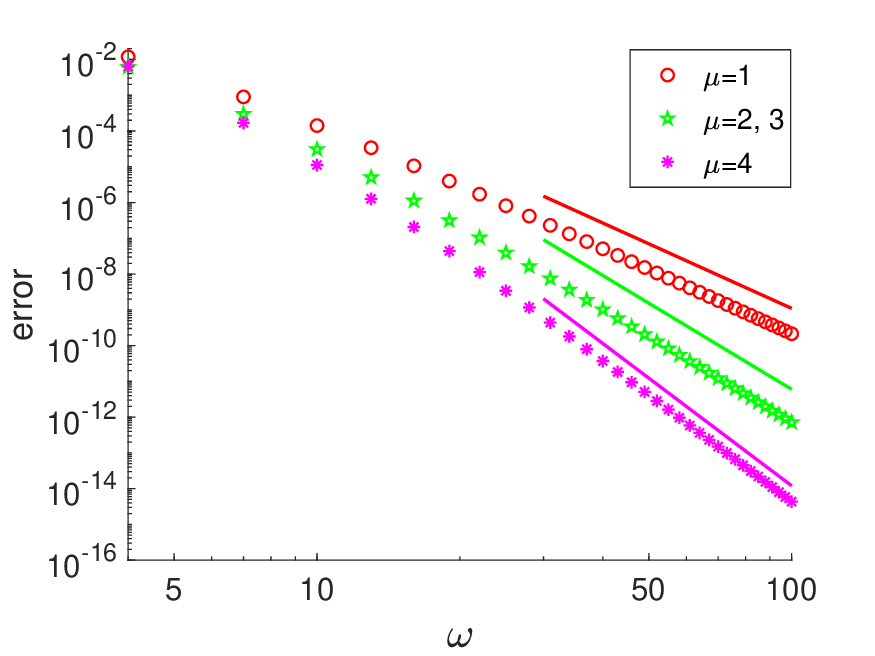}~~
\includegraphics[width=7.4cm,height=6.2cm]{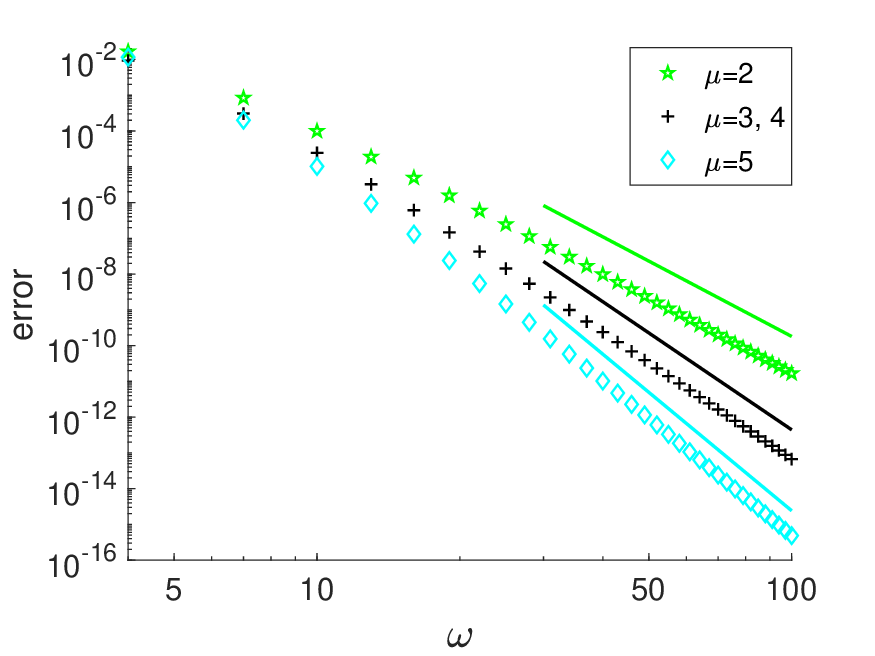}\\
\caption{Absolute errors of $(\mathcal{Q}_{2n,\mu}^{\mathrm{HI}}f)(\omega)$ with $n=1$ as a function of $\omega$ for $\nu=1$ (left) and $\nu=2$ (right).  Here $f(x)= e^{-x}$ (top) and $f(x)=1/(1+x)^2$ (bottom), we choose different values of $\mu$ and the solid lines indicate the predicted rates $\mathcal{O}(\omega^{-4n-\mu-1})$ when $\mu-\nu$ is even and $\mathcal{O}(\omega^{-4n-\mu-2})$ when $\mu-\nu$ is odd.} \label{fig:GRinteger1}
\end{figure}

\begin{figure}
\centering
\includegraphics[width=7.4cm,height=6.2cm]{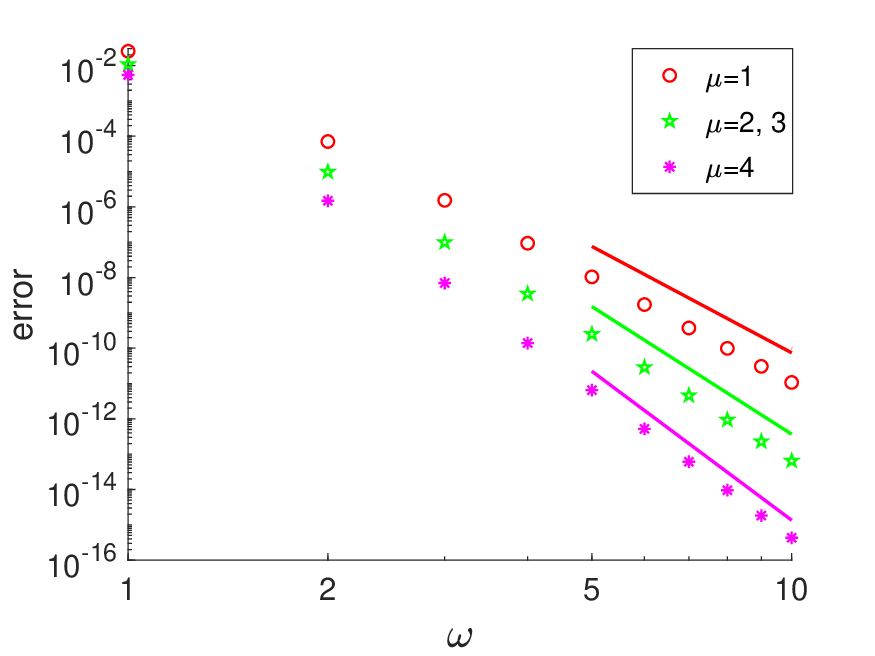}~~
\includegraphics[width=7.4cm,height=6.2cm]{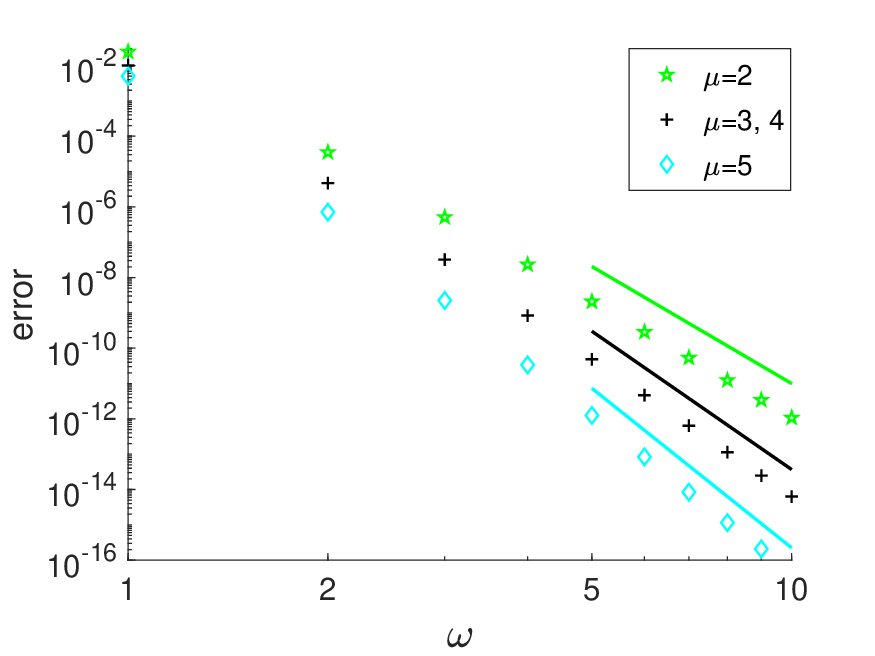}\\
\includegraphics[width=7.4cm,height=6.2cm]{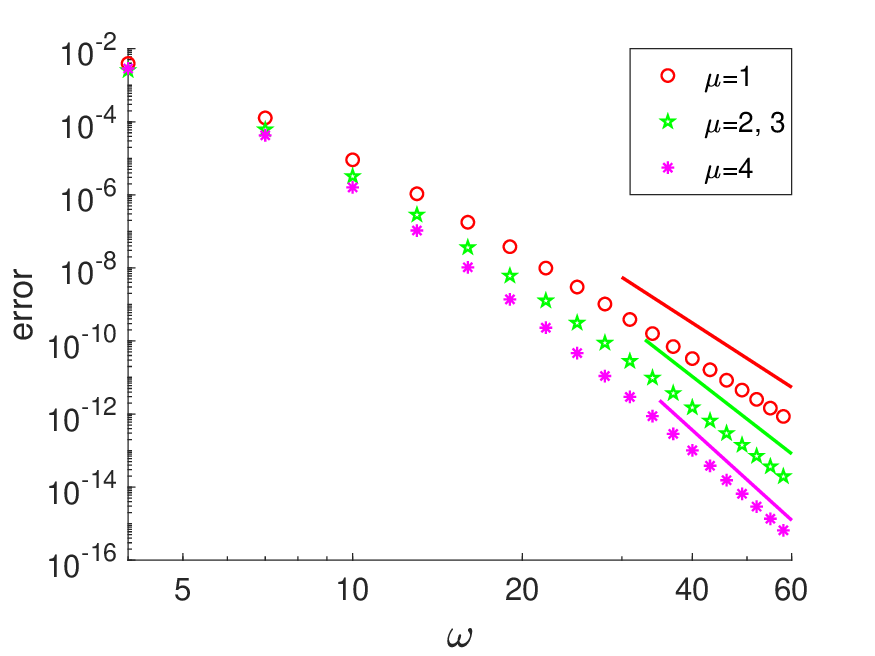}~~
\includegraphics[width=7.4cm,height=6.2cm]{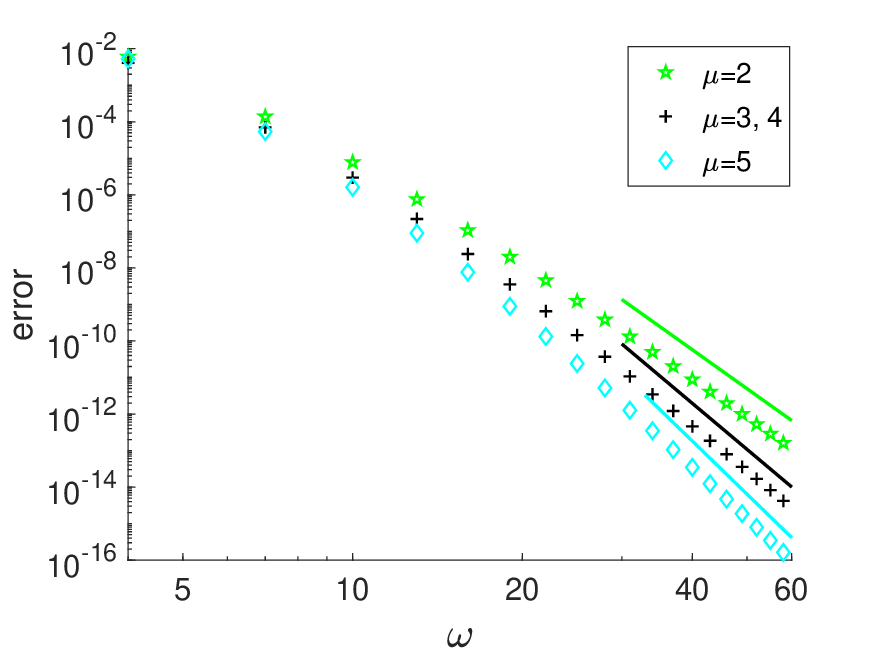}\\
\caption{Absolute errors of $(\mathcal{Q}_{2n,\mu}^{\mathrm{HI}}f)(\omega)$ with $n=2$ as a function of $\omega$ for $\nu=1$ (left) and $\nu=2$ (right).  Here $f(x)= e^{-x}$ (top) and $f(x)=1/(1+x)^2$ (bottom), we choose different values of $\mu$ and the solid lines indicate the predicted rates $\mathcal{O}(\omega^{-4n-\mu-1})$ when $\mu-\nu$ is even and $\mathcal{O}(\omega^{-4n-\mu-2})$ when $\mu-\nu$ is odd.} \label{fig:GRinteger2}
\end{figure}

We demonstrate the performance of the proposed rule $(\mathcal{Q}_{2n,\mu}^{\mathrm{HI}}f)(\omega)$ in Theorem \ref{them:HankelInteger}. In Figures \ref{fig:GRinteger1} and \ref{fig:GRinteger2} we plot the absolute errors of $(\mathcal{Q}_{2n,\mu}^{\mathrm{HI}}f)(\omega)$ as a function of $\omega$ for $n=1$ and $n=2$, respectively. For each $\nu$, we consider several different values of $\mu$. We see that the errors of $(\mathcal{Q}_{2n,\mu}^{\mathrm{HI}}f)(\omega)$ decay at the rate $\mathcal{O}(\omega^{-4n-\mu-1})$ when $\mu-\nu$ is even and the rate $\mathcal{O}(\omega^{-4n-\mu-2})$ when $\mu-\nu$ is odd, which are consistent with the error estimate in Theorem \ref{them:GGRQ}.

\begin{figure}
\centering
\includegraphics[width=7.4cm,height=6.2cm]{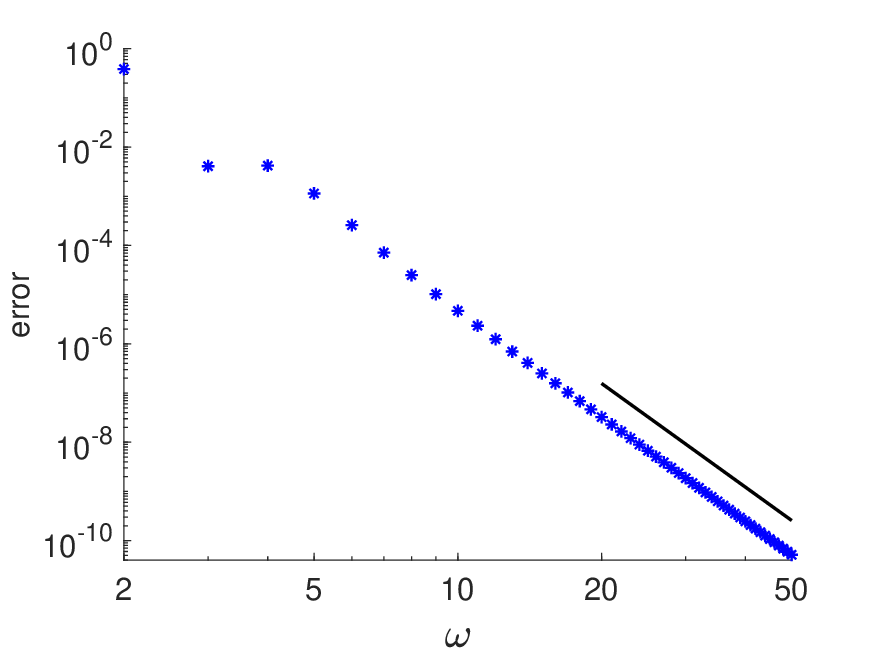}~~~
\includegraphics[width=7.4cm,height=6.2cm]{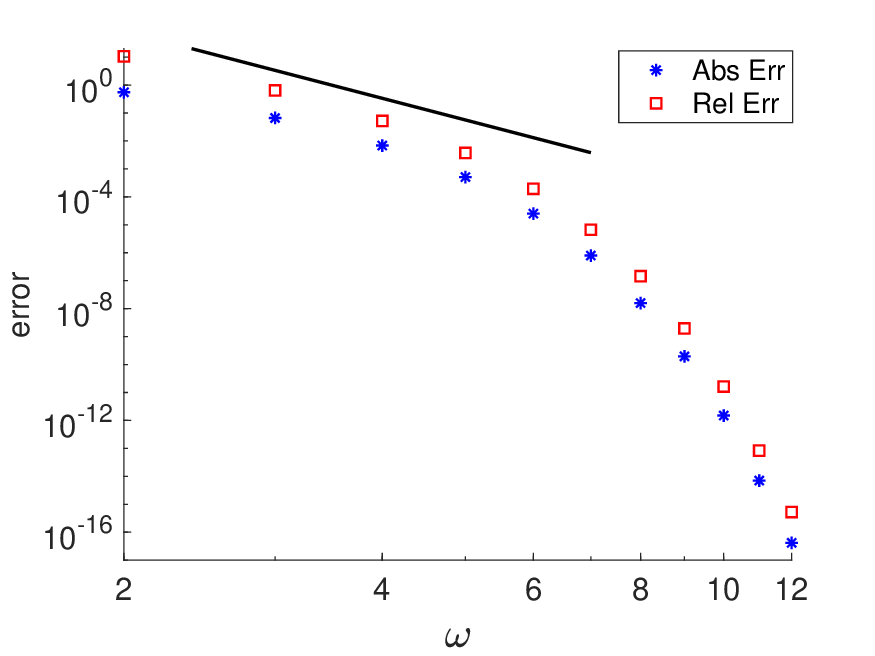}
\caption{Absolute errors of $(\mathcal{Q}_{2n,\mu}^{\mathrm{HI}}f)(\omega)$ as a function of $\omega$ for $n=1$, $\nu=2$ (left) and absolute and relative errors of $(\mathcal{Q}_{2n,\mu}^{\mathrm{HI}}f)(\omega)$ for $n=1$, $\nu=3$ (right). Here $f(x)=e^{-x^2}$, $\mu=\nu$ and the solid lines indicate the predicted rates $\mathcal{O}(\omega^{-4n-\mu-1})$.}
\label{fig:specialcase}
\end{figure}

We mention an interesting superconvergence phenomenon of the rule $(\mathcal{Q}_{2n,\mu}^{\mathrm{HI}}f)(\omega)$ in certain special situations. In Figure \ref{fig:specialcase} we display the absolute error of the rule as a function of $\omega$ for $f(x)=e^{-x^2}$ and $\nu=2$ and the absolute and relative errors
of the rule for $\nu=3$ and we choose $\mu=\nu$ in our calculations. We see that $(\mathcal{Q}_{2n,\mu}^{\mathrm{HI}}f)(\omega)$ converges at the predicted rate when $\nu$ is even, but at a much faster rate when $\nu$ is odd. In fact, by \cite[Equation (6.618,1)]{Gradshteyn2007} we know for $\nu>-1$ that
\begin{equation}
\int_{0}^{\infty} e^{-x^2} J_{\nu}(\omega x) \ud x =
\frac{\sqrt{\pi}}{2} \exp\left(-\frac{\omega^2}{8}\right) I_{\frac{\nu}{2}}\left(\frac{\omega^2}{8}\right), \nonumber
\end{equation}
where $I_{\nu}(z)$ is the modified Bessel function of the first kind. When $\nu$ is odd, using the above closed form and noting that $f(x)$ is even, we find after some elementary calculations that the convergence rate of $(\mathcal{Q}_{2n,\mu}^{\mathrm{HI}}f)(\omega)$ is $O(e^{-\omega^2}/\omega)$ as $\omega\rightarrow\infty$, which explains the superconvergence phenomenon displayed in Figure \ref{fig:specialcase}. However, we still do not know how to give a precise condition on $f(x)$ and $\nu$ such that the proposed rule has the superconvergence phenomenon.

\begin{figure}
\centering
\includegraphics[width=0.45\textwidth]{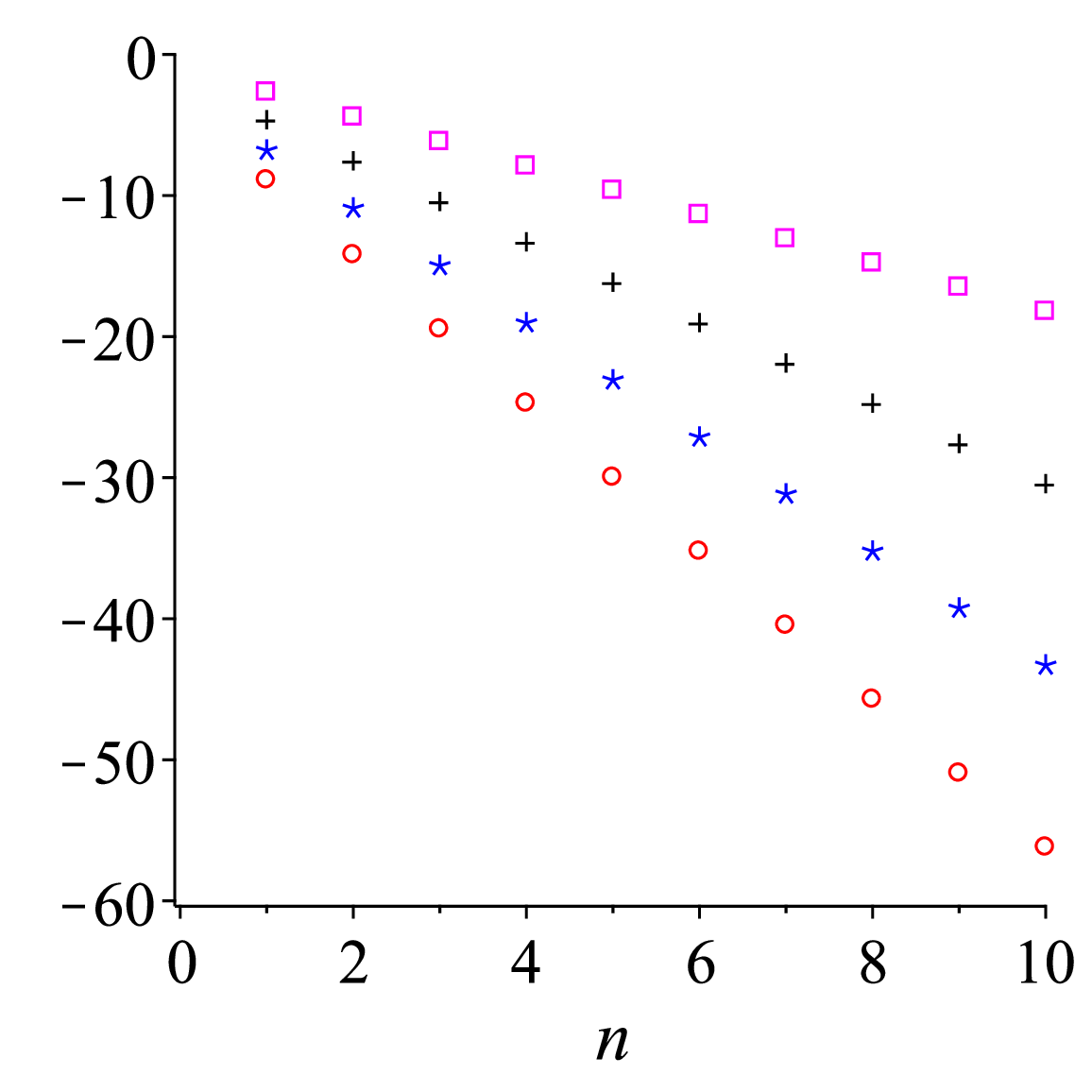}~~~~~~
\includegraphics[width=0.45\textwidth]{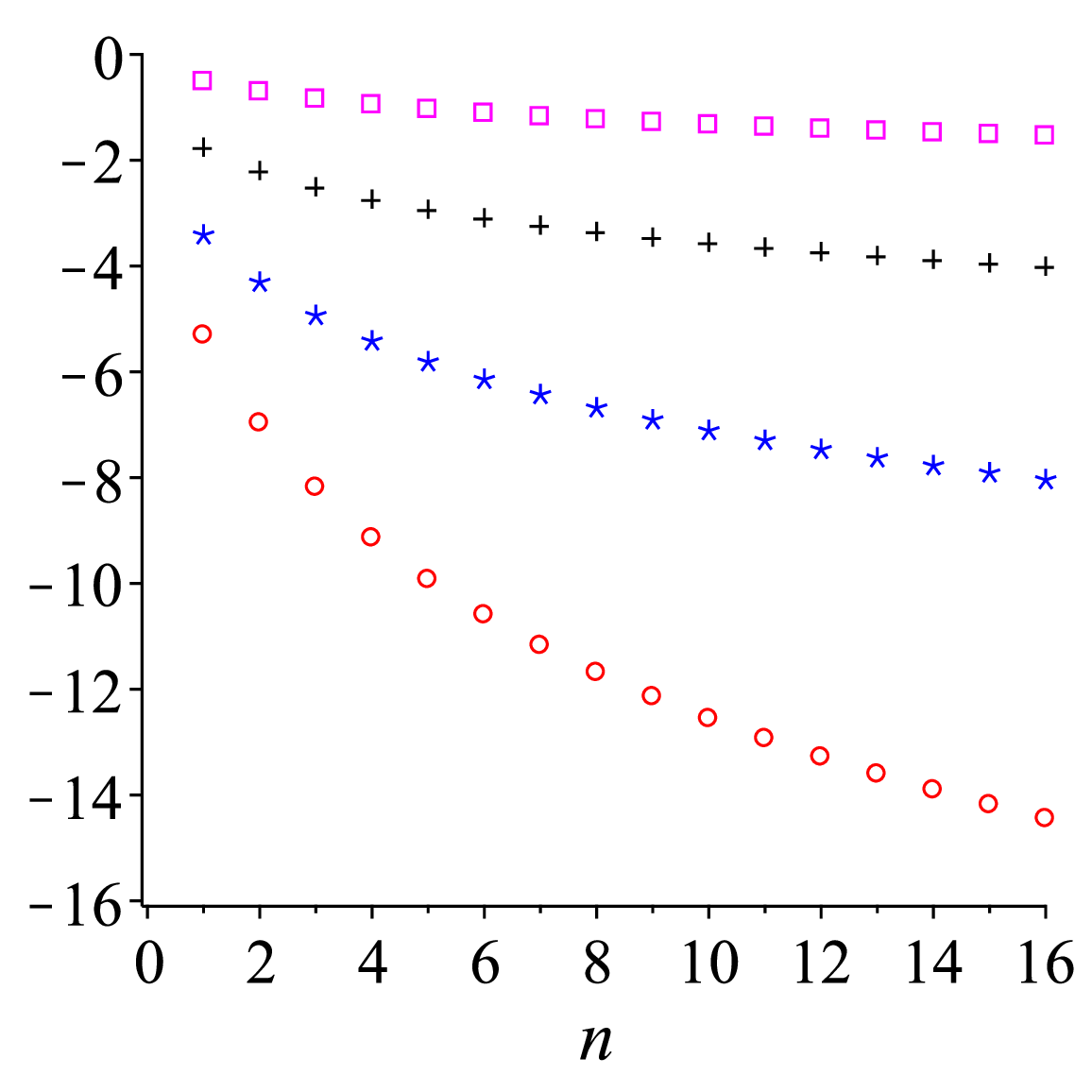}
\caption{The logarithm of the absolute error of $(\mathcal{Q}_{2n,\mu}^{\mathrm{HI}}f)(\omega)$ as a function of $n$ for $\nu=2$ and $\mu=2$ and $f(x)= e^{-x}$ (left) and $f(x)= 1/(1+x)^2$ (right). Here $\omega=2$ ($\Box$), $\omega=4$ ($+$), $\omega=8$ ($\star$), $\omega=16$ ($\circ$).}
\label{fig:nlarge}
\end{figure}

Finally, we demonstrate in Figure \ref{fig:nlarge} the accuracy of the rule $(\mathcal{Q}_{2n,\mu}^{\mathrm{HI}}f)(\omega)$ as a function of $n$. For each fixed $\omega$, we see that the accuracy of $(\mathcal{Q}_{2n,\mu}^{\mathrm{HI}}f)(\omega)$ improves as $n$ increases and the larger $\omega$, the faster the accuracy improved.

\section{Applications}\label{sec:Appl}
In this section, we present two applications of the generalized Gauss-Radau quadrature rule in Theorem \ref{them:HankelInteger}.
\subsection{Oscillatory Hilbert transform}
Consider the following oscillatory Hilbert transform
\begin{equation}
(H_{\nu}f)(\omega,\tau) := \dashint_{0}^{\infty} \frac{f(x)}{x-\tau} J_{\nu}(\omega x) \ud x,
\end{equation}
where $\tau>0$, $\nu\in\mathbb{N}_0$ and the bar indicates the Cauchy principal value. To construct efficient methods for computing such transform, the main difficulties are that the integrand is oscillatory and has a singularity of Cauchy-type.

In the following we present a method for computing this transform. By the singularity subtraction technique, we have
\begin{equation}\label{eq:HilbSub}
(H_{\nu}f)(\omega,\tau) := \int_{0}^{\infty} \frac{f(x)-f(\tau)}{x-\tau} J_{\nu}(\omega x) \ud x + f(\tau) \dashint_{0}^{\infty} \frac{J_{\nu}(\omega x)}{x-\tau}  \ud x.
\end{equation}
For the last integral, it can be expressed by the Struve or Meijer G-functions (see, e.g., \cite{Wang2013,Xu2016}) and their evaluations can be performed by most modern software, e.g., Maple and Matlab. As for the first integral on the right-hand side of \eqref{eq:HilbSub}, it can be evaluated directly by using the generalized Gauss-Radau quadrature rule in Theorem \ref{them:HankelInteger}. When the singularity $\tau$ is not close to zero, an obvious advantage of using this quadrature rule is that it can avoid the loss of accuracy due to cancellation since its nodes all lie on the imaginary axis.

\begin{figure}
\centering
\includegraphics[width=7.4cm,height=6.2cm]{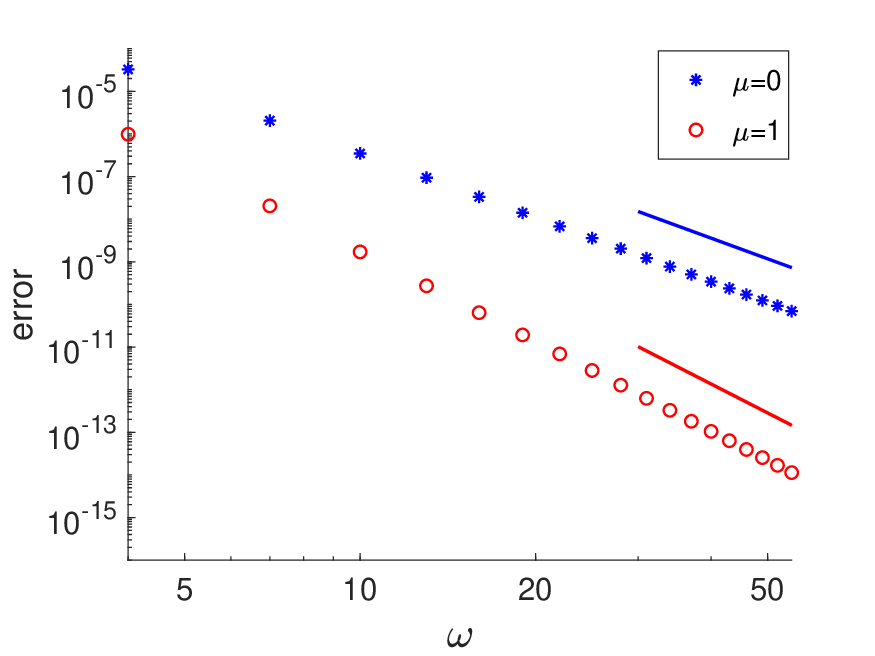}~~~
\includegraphics[width=7.4cm,height=6.2cm]{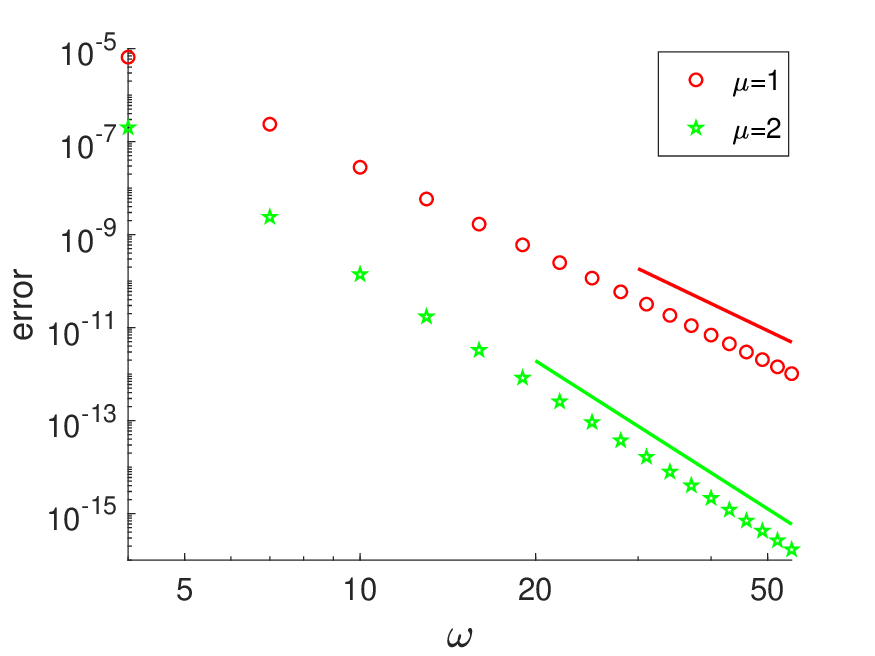}\\
\includegraphics[width=7.4cm,height=6.2cm]{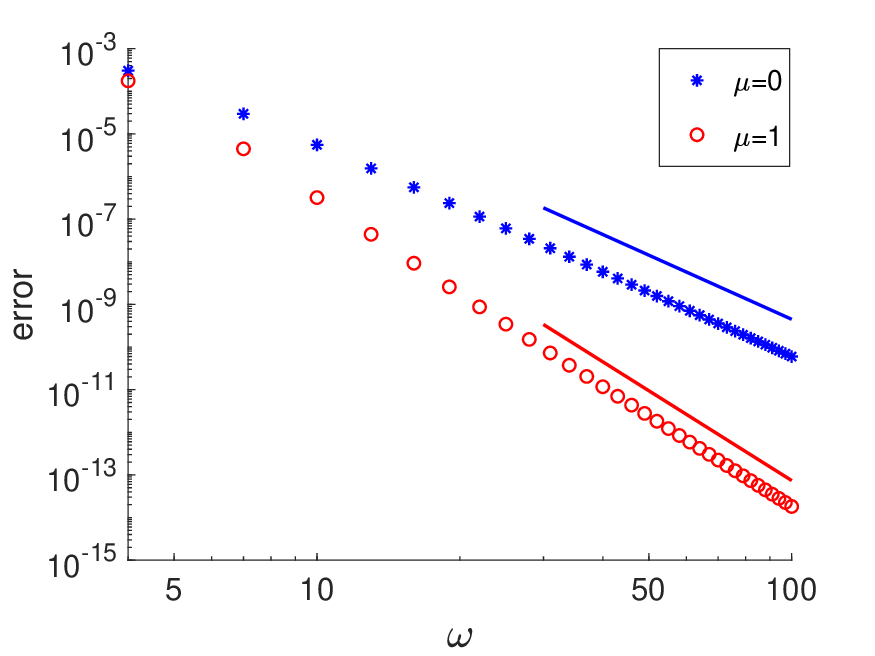}~~~
\includegraphics[width=7.4cm,height=6.2cm]{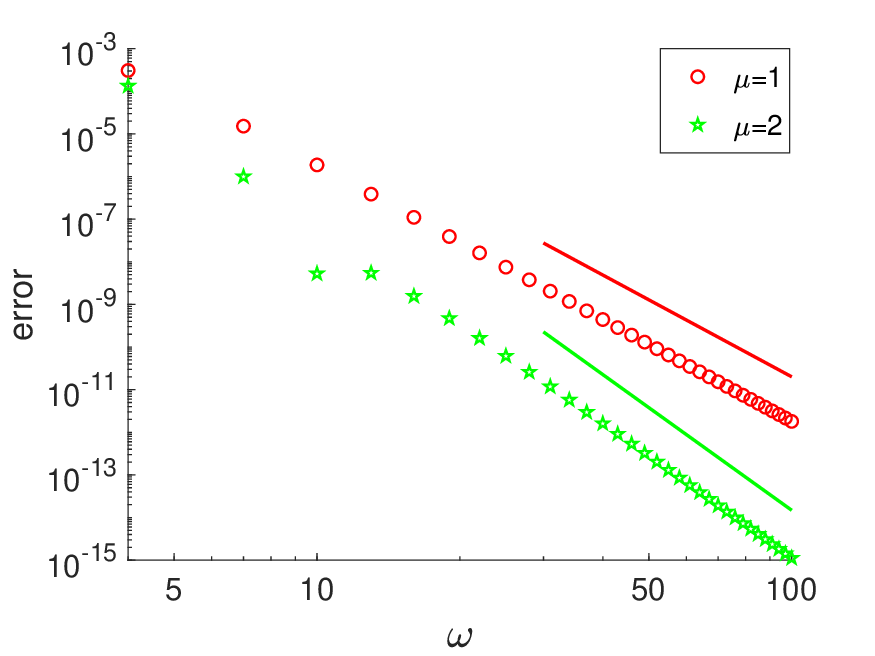}
\caption{Absolute errors as a function of $\omega$ for $\nu=0$ (left) and $\nu=1$ (right). Top row shows $f(x)= e^{-x}$ and $\tau=5$ and bottom row shows $f(x)=1/(1+(1+x)^2)$ and $\tau=1$ and the solid lines indicate the predicted rates.}
\label{fig:Hilbert}
\end{figure}

To show the performance of the above proposed method, we consider the test functions $f(x)=e^{-x}$ and $f(x)=1/(1+(1+x)^2)$, the orders $\nu=0,1$ and the singularities $\tau=1,5$. We evaluate the first integral on the right-hand side of \eqref{eq:HilbSub} by the rule $(\mathcal{Q}_{2n,\mu}^{\mathrm{HI}}f)(\omega)$ with $n=1$ and the second integral by \cite[Corollary 4.2]{Wang2013}\footnote{We point out that \cite[Equation (4.8)]{Wang2013} is correct for $\nu=0$, but is wrong for $\nu=1$ since the integral in \cite[Equation (4.2)]{Wang2013} will involve a nonintegrable singularity when setting $f(t)=1$. However, by setting $f(t)=t$ in \cite[Equation (4.2)]{Wang2013}, it is not difficult to derived the correct result given here.}
\begin{equation}\label{eq:ClosedForm}
\left.
\begin{array}{lll}
{\displaystyle \dashint_{0}^{\infty} \frac{J_{0}(\omega x)}{x-\tau} \ud x = -\frac{\pi}{2} \big[ \mathbf{H}_0(\omega\tau) + Y_0(\omega\tau) \big] },  \\[3ex]
{\displaystyle \dashint_{0}^{\infty} \frac{J_{1}(\omega x)}{x-\tau} \ud x =  \frac{\pi}{2} \big[ \mathbf{H}_{-1}(\omega\tau) - Y_1(\omega\tau) \big] - \frac{1}{\omega\tau} },
\end{array}
\right.
\end{equation}
where $\mathbf{H}_{\nu}(z)$ is the Struve function and $Y_{\nu}(z)$ is the Bessel function of the second kind. For each $\nu$, we choose $\mu=\nu$ and $\mu=\nu+1$. Here we ignore the errors of computing the integrals \eqref{eq:ClosedForm} and the error of the proposed method comes only from the rule $(\mathcal{Q}_{2n,\mu}^{\mathrm{HI}}f)(\omega)$ and thus it decays at the rate $\mathcal{O}(\omega^{-4n-\mu-1})$ when $\mu-\nu$ is even and the rate $\mathcal{O}(\omega^{-4n-\mu-2})$ when $\mu-\nu$ is odd. Figure \ref{fig:Hilbert} plots the results and we see that the proposed method converge at expected rates.

\subsection{Electromagnetic field configurations}
In electromagnetic geophysics, the measurements of the electromagnetic field components is achieved by placing magnetic dipoles above the surface which consist of transmitter and receiver loops that can be horizontal co-planar and perpendicular \cite{Denich2023,Ward1987}. Suppose that the local subsurface below the instrument composed by horizontal and homogeneous layers and let $N$ denote the number of layers. The components of the magnetic field for horizontal co-planar and perpendicular configurations are given respectively by
\begin{align}\label{eq:Horizon}
H_{z}^{(N)} &= \frac{m}{4\pi} \int_{0}^{\infty} (1 + \Phi_0(\lambda)) \lambda^2 J_{0}(\omega \lambda) \ud \lambda, \nonumber \\[2ex]
H_{\rho}^{(N)} &= \frac{m}{4\pi} \int_{0}^{\infty} (1 - \Phi_0(\lambda)) \lambda^2 J_{1}(\omega \lambda) \ud \lambda, \nonumber
\end{align}
where $m$ is the transmitter's magnetic moment, $\omega$ is the distance between the transmitter and the receiver and $\Phi_0(\lambda)$ is the reflection term defined recursively by
\begin{align}
&\Phi_N(\lambda) = 0, \quad \Phi_j(\lambda) = \frac{\Phi_{j+1}(\lambda)+\Psi_{j+1}(\lambda)}{\Phi_{j+1}(\lambda)\Psi_{j+1}(\lambda)+1}e^{-2u_j(\lambda)h_j}, \quad j=N-1,\ldots,1,  \nonumber \\
&\Phi_0(\lambda) = \frac{\Phi_1(\lambda)+\Psi_1(\lambda)}{\Phi_1(\lambda)\Psi_1(\lambda)+1},  \nonumber
\end{align}
and $\Psi_j(\lambda)$ and $u_j(\lambda)$ are given by
\begin{align}
\Psi_j(\lambda) = \frac{u_{j-1}(\lambda) - u_j(\lambda)}{u_{j-1}(\lambda) + u_j(\lambda)}, \quad  u_j(\lambda)=\sqrt{\lambda^2-k_j^2}, \quad j=1,\ldots,N, \nonumber
\end{align}
and $u_0(\lambda)=\lambda$ and $k_{j}=\sqrt{-\mathrm{i}\omega_{0}\mu_{0}\sigma_{j}}$ for $j=1,\ldots,N$. Here $\mu_{0}=4\pi\times10^{-7}~\textrm{H/m}$ is the magnetic permeability of vacuum, $\omega_{0} = 2 \pi f_{0}$ is the angular frequency and $f_{0}$ is the transmitter's frequency, $\sigma_{j}$ and $h_{j}$ represent conductivity and thickness of the $j$-th layer.

\begin{figure}
\centering
\includegraphics[width=7.4cm,height=6.2cm]{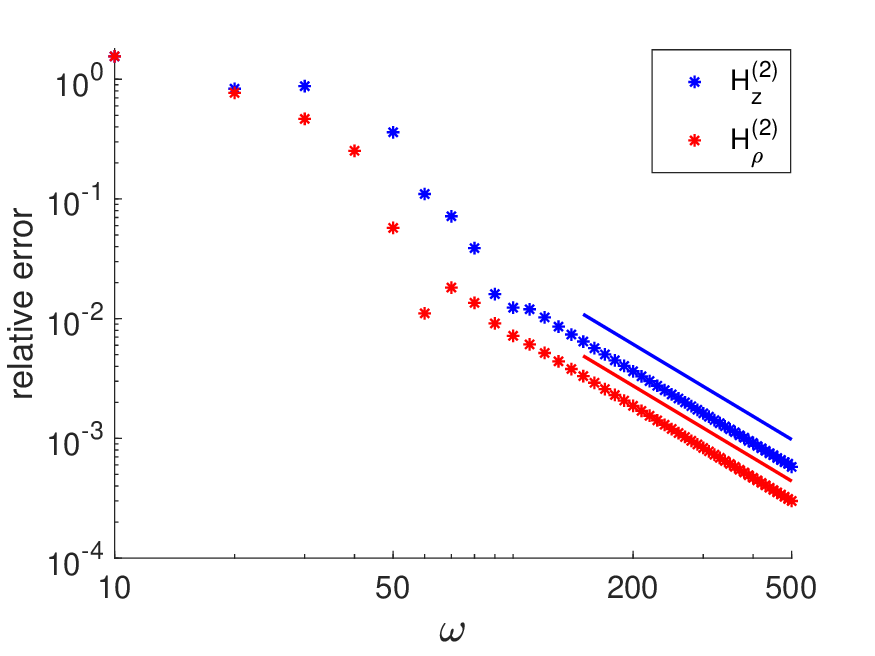}~~~
\includegraphics[width=7.4cm,height=6.2cm]{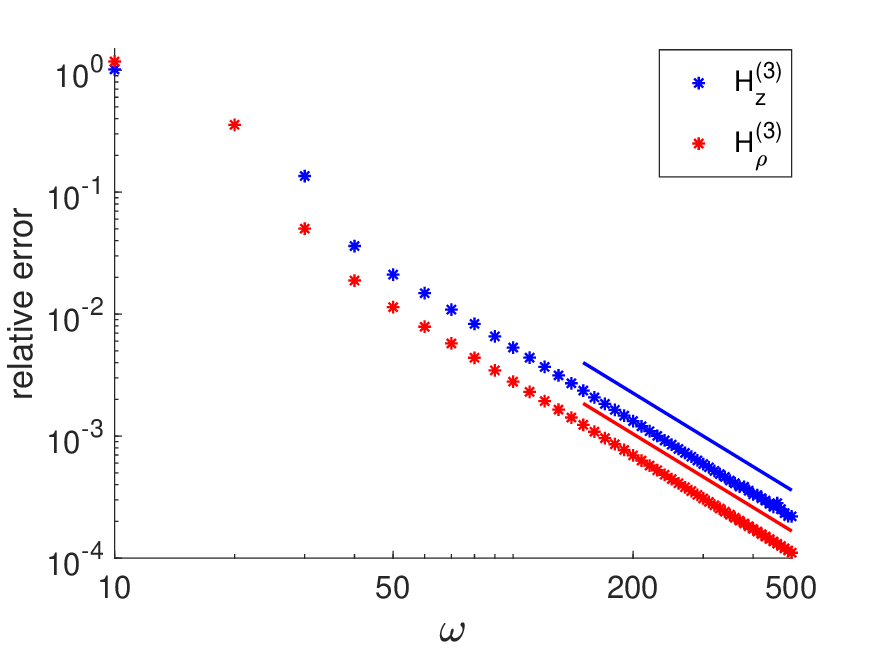}\\
\caption{Relative errors of $(\mathcal{Q}_{2n,\mu}^{\mathrm{HI}}f)(\omega)$ as a function of $\omega$ for $H_{z}^{(N)}$ and $H_{\rho}^{(N)}$ with $N=2$ (left) and $N=3$ (right).
Here we choose $n=\mu=1$, $m=1$, $f_{0}=1000$, and
$\sigma_{1}=50$, $\sigma_{2}=4.9$, $h_{1}=3$ for $N=2$ and $\sigma_{1}=76.9$, $\sigma_{2}=32.3$, $\sigma_{3}=50$, $h_{1}=2.5$, $h_{2}=0.5$ for $N=3$. The solid lines indicate the rates $\mathcal{O}(\omega^{-2})$.}
\label{fig:layer}
\end{figure}

Note that $H_{z}^{(N)}$ and $H_{\rho}^{(N)}$ are Hankel transforms of orders zero and one, respectively. We consider the problem of computing them with $N=2,3$ by using the rule $(\mathcal{Q}_{2n,\mu}^{\mathrm{HI}}f)(\omega)$ in \eqref{eq:GGaussRadau}. By Theorem \ref{them:GGRQ} we know that the absolute errors of $(\mathcal{Q}_{2n,\mu}^{\mathrm{HI}}f)(\omega)$ are $\mathcal{O}(\omega^{-7})$ for $H_{z}^{(N)}$ and $\mathcal{O}(\omega^{-6})$ for $H_{\rho}^{(N)}$ as $\omega\rightarrow\infty$. Moreover, note that $H_{z}^{(N)}=\mathcal{O}(\omega^{-5})$ and $H_{\rho}^{(N)}=\mathcal{O}(\omega^{-4})$ as $\omega\rightarrow\infty$, and therefore we can expect that the relative errors of the rule $(\mathcal{Q}_{2n,\mu}^{\mathrm{HI}}f)(\omega)$ decay at the rate $\mathcal{O}(\omega^{-2})$ as $\omega\rightarrow\infty$. In Figure \ref{fig:layer} we plot the relative errors of the rule $(\mathcal{Q}_{2n,\mu}^{\mathrm{HI}}f)(\omega)$ with $n=\mu=1$ as a function of $\omega$ for $N=2,3$. We see that numerical results are consistent with our prediction.

\section{Conclusion}\label{sec:concluding}
In this paper, we have studied the construction of complex generalized Gauss-Radau quadrature rules for Hankel transform of integer order. We have shown that, by adding certain function and derivative values at the left endpoint, complex generalized Gauss-Radau quadrature rules can be constructed with guaranteed existence. Motivated by this finding, we further studied the existence of the polynomials that are orthogonal with respect to the oscillatory weight function $x^{\mu}J_{\nu}(x)$ on $(0,\infty)$. When $\mu-\nu$ is a nonnegative integer, we proved that the polynomials always exist for all degrees if $\mu-\nu$ is even, but exist only for even degrees if $\mu-\nu$ is odd.

Gaussian quadrature rules for highly oscillatory integrals have received increasing attention in recent years due to the fact that they can achieve optimal asymptotic order \cite{Asheim2013,Asheim2014,Huybrechs2019,Wu2024}. A key issue of such rules is that the existence of the orthogonal polynomials with respect to the sign-changing weight function cannot be guaranteed and the proofs are generally quite challenging (see, e.g., \cite{Celsus2022,Deano2016,Huybrechs2019}). The findings of the present study give a sequence of such polynomials with guaranteed existence, which provide a theoretical rationale for the developed generalized Gauss-Radau quadrature rules. 

\section*{Acknowledgements}
We are grateful to the reviewers for their helpful comments to improve the presentation of the paper. The first author would like to thank Daan Huybrechs for his valuable feedback on polynomials that are orthogonal with respect to Bessel functions and for helpful discussions on generalized Gaussian quadrature rules for oscillatory integrals during the workshop "Singular and Oscillatory Integration: Advances and Applications" held at University College London from June 24th to 26th, 2024. This work was supported by the National Natural Science Foundation of China under grant number 12371367.

\end{document}